\documentclass[a4paper]{article}

\usepackage{graphicx}
\usepackage{color}
\usepackage{stmaryrd}

\usepackage{tikz}
\usepackage{tikz-cd}
\usetikzlibrary{positioning,calc}
\usepackage{ifthen}

\newlength{\dsunit} % Scaling unit on which all other lengths are
\newlength{\nodedist} % Distance between nodes, as used by relative
\newlength{\dsmin} % Minimum height and width for boxes
\newcommand{\dirout}[4]{ 
\node[outer] [#4=2\dsunit of #2] {$\scriptstyle
\ifthenelse{\equal{#3}{}}{}{\vphantom{g\ell}}
#3$}
edge[-] (#1);
}

\newcommand{\outabovex}[3]{\dirout{#1}{#2}{#3}{above}}
\newcommand{\outrightx}[3]{\dirout{#1}{#2}{#3}{right}}
\newcommand{\outbelowx}[3]{\dirout{#1}{#2}{#3}{below}}
\newcommand{\outleftx}[3]{\dirout{#1}{#2}{#3}{left}}
\newcommand{\outabove}[2]{\outabovex{#1}{#1}{#2}}
\newcommand{\outright}[2]{\outrightx{#1}{#1}{#2}}
\newcommand{\outbelow}[2]{\outbelowx{#1}{#1}{#2}}
\newcommand{\outleft}[2]{\outleftx{#1}{#1}{#2}}

\newcommand{\inedge}[3]{
\draw[-] (#1) to node{$\scriptstyle #3$} (#2);
}
\newcommand{\inedges}[3]{
\draw[-] (#1) to node[swap]{$\scriptstyle #3$} (#2);
}

\newcommand{\throughvert}[3]{
\node[ghost] (#1up) [above=2\dsunit of #1] {$\scriptstyle #3$};
\node[ghost] (#2dn) [below=2\dsunit of #2] {$\scriptstyle #3$};
\draw[-] (#1up) to (#2dn);
}

\newcommand{\throughhoriz}[3]{
\node[ghost] (#1up) [left=2\dsunit of #1] {$\scriptstyle #3$};
\node[ghost] (#2dn) [right=2\dsunit of #2] {$\scriptstyle #3$};
\draw[-] (#1up) to (#2dn);
}

\newenvironment{doublestring}{%
\begin{array}{c}
\begin{tikzpicture}
[auto,
x=10\dsunit,y=10\dsunit,
node distance=\nodedist, % manual p.54
ghostdebug/.style={rectangle,fill=yellow,minimum width=\dsmin,minimum
height=\dsmin}, % for debugging
ghost/.style={rectangle,minimum width=\dsmin,minimum height=\dsmin},
outer/.style={rectangle,inner sep=1pt}, 
basic/.style={rectangle,inner sep=1pt,draw,minimum width=\dsmin,minimum height=\dsmin},
w2/.style={rectangle,inner sep=1pt,draw,minimum width=2*\dsmin+\nodedist,minimum height=\dsmin},
w3/.style={rectangle,inner sep=1pt,draw,minimum width=3*\dsmin+2*\nodedist,minimum height=\dsmin},
w4/.style={rectangle,inner sep=1pt,draw,minimum width=4*\dsmin+3*\nodedist,minimum height=\dsmin},
w6/.style={rectangle,inner sep=1pt,draw,minimum width=6*\dsmin+5*\nodedist,minimum height=\dsmin},
h2/.style={rectangle,inner sep=1pt,draw,minimum width=\dsmin,minimum height=2*\dsmin+\nodedist},
h3/.style={rectangle,inner sep=1pt,draw,minimum width=\dsmin,minimum height=3*\dsmin+2*\nodedist},
h4/.style={rectangle,inner sep=1pt,draw,minimum width=\dsmin,minimum height=4*\dsmin+3*\nodedist}
]}
{\end{tikzpicture}\end{array}}

\usepackage[numbers]{natbib}
\usepackage{natbibspacing}
\renewcommand{\bibfont}{\small}

\definecolor{mylinkcolor}{rgb}{0.05,0.05,0.4}
\usepackage[linktocpage]{hyperref}
\hypersetup{colorlinks,
linkcolor=mylinkcolor,
citecolor=mylinkcolor,
urlcolor=mylinkcolor}

\usepackage{latexsym}
\usepackage{amsmath,amssymb}
\newcommand{\demphcolor}{\color[rgb]{0,0,0}}
\newcommand{\demph}[1]{\textbf{\textup{\demphcolor#1}}}

\newcommand{\fcat}[1]{\mathbf{#1}}
\newcommand{\scat}[1]{\mathbf{#1}}
\newcommand{\ovln}[1]{\overline{#1}}
\newcommand{\twid}[1]{\widetilde{#1}}
\renewcommand{\epsilon}{\varepsilon}
\newcommand{\iso}{\cong}
\newcommand{\eqv}{\simeq}
\newcommand{\of}{\mathbin{\circ}}

\newcommand{\ladj}{\dashv}
\newcommand{\from}{\colon}

\newcommand{\textimplies}{\Rightarrow}

\newcommand{\Z}{\mathbb{Z}}

\newcommand{\op}{\mathrm{op}}

\newcommand{\Cat}{\fcat{Cat}}

\newcommand{\oppairu}{\rightleftarrows}

\newcommand{\swnt}{\rotatebox{45}{$\Leftarrow$}\!}

\newcommand{\toby}[1]{\xrightarrow{#1}}
\newcommand{\otby}[1]{\xleftarrow{#1}}

\newcommand{\oppair}[4]{%
\begin{tikzcd}[cramped, sep=small, ampersand replacement=\&] 
#1
\ar[r, shift left, "#3"] \& 
#2
\ar[l, shift left, "#4"] 
\end{tikzcd}}

\def\today{\number\day\space \ifcase\month\or
  January\or February\or March\or April\or May\or June\or
  July\or August\or September\or October\or November\or December\fi
  \space\number\year}

\renewcommand{\phi}{\varphi}
\def\mathcal{\mathscr}

\newcommand{\ot}{\leftarrow}

\newcommand{\A}{\scat{A}}
\newcommand{\B}{\scat{B}}
\newcommand{\C}{\scat{C}}
\newcommand{\Ho}{H} %{\textup{H}}
\newcommand{\Ve}{V} %{\textup{V}}
\newcommand{\vertmap}[3]{\begin{tikzcd}[cramped, sep=small]
#1 \ar[d, "#3"]\\ #2\end{tikzcd}}

\usepackage[amsmath,thmmarks]{ntheorem}

\newtheorem{thm}{Theorem}[section]

\newtheorem{lemma}[thm]{Lemma}
\newtheorem{cor}[thm]{Corollary}

\newtheorem{principle}[thm]{Principle}

\newtheorem*{thmA}{Theorem~\ref{thm:A}}
\newtheorem*{thmB}{Theorem~\ref{thm:B}}
\newtheorem*{thmC}{Theorem~\ref{thm:C}}
\newtheorem*{thmD}{Theorem~\ref{thm:D}}

\theorembodyfont{\normalfont}

\newtheorem{defn}[thm]{Definition}
\newtheorem{example}[thm]{Example}

\newtheorem{remark}[thm]{Remark}

\theoremstyle{nonumberplain}
\theoremsymbol{\ensuremath{\square}}
\qedsymbol{\ensuremath{\square}}

\newtheorem{proof}{Proof}

\newcommand{\theoremtobeproved}{}
\newtheorem{pfoftheorem}{Proof of \theoremtobeproved}
\newenvironment{pfof}[1]
{
\renewcommand{\theoremtobeproved}{#1}
\begin{pfoftheorem}
}
{\end{pfoftheorem}}

\allowdisplaybreaks

\title{Equivalence via surjections}
\author{Tom Leinster
\thanks{Maxwell Institute and School of Mathematics, University of
Edinburgh, Scotland; Tom.Leinster@ed.ac.uk}}
\date{}

\begin{document}

\sloppy
\maketitle

\begin{abstract}
Many types of categorical structure obey the following principle: the
natural notion of equivalence is generated, as an equivalence relation, by
identifying $\A$ with $\B$ when there exists a strictly
structure-preserving map $\A \to \B$ that is genuinely (not just
essentially) surjective in each dimension and faithful in the top
dimension. We prove this principle for four types of structure: categories,
monoidal categories, bicategories and double categories. The last of these
theorems suggests that the right notion of equivalence between double
categories is Campbell's gregarious double equivalence, a conclusion also
reached for different reasons in recent work of Moser, Sarazola and
Verdugo.
\end{abstract}

\section{Introduction}
\label{sec:intro}

This note was inspired by the invited lecture of Maru Sarazola at Category
Theory 2025, which led with the question `What is the right notion of
equivalence between double categories?'~\cite{SaraDCE}. Sarazola's answer
was that it is Alexander Campbell's gregarious double
equivalence~\cite{Camp}, and she gave a model-categorical justification
representing forthcoming joint work with Lyne Moser and Paula
Verdugo~\cite{MSV}, adding to the argument of Campbell himself.

Here I give a different reason why gregarious double equivalence is the
right answer. (At least, it is right subject to a subtlety concerning
strictness.) The argument is of interest beyond the context of double
categories, and represents one instance of the following principle.

\begin{principle}
\label{princ:main}
The natural notion of equivalence between categorical structures is
generated by identifying structures $\A$ and $\B$ when there exists a strict
surjective equivalence $\A \to \B$.
\end{principle}

The meaning of this principle is best explained by example. 

The simplest instance, and the first of the four theorems we prove,
concerns categories themselves. There, \emph{surjective equivalence} means
a functor that is genuinely surjective on objects (not just essentially
surjective), full, and faithful. The word \emph{strict} does nothing in
this simplest instance.

Let $\sim$ be the equivalence relation on the class of categories generated
by declaring $\A \sim \B$ if there exists a surjective equivalence $\A \to
\B$. Thus, $\A \sim \B$ if and only if there exists a zigzag of surjective
equivalences
\[
\A = \A_0 \ot \A_1 \to \A_2 \ot \ \cdots \to \A_n = \B.
\]
The class of surjective equivalences is stable under pullback, so a
simpler description is this: $\A \sim \B$ if and only if there exists a
span
\[
\A \ot \C \to \B
\]
of surjective equivalences. Principle~\ref{princ:main} correctly suggests
that $\sim$ is the standard notion of equivalence of categories:

\begin{thmA}
Two categories $\A$ and $\B$ are equivalent if and only if there exist a
category $\C$ and surjective equivalences $\A \ot \C \to \B$.
\end{thmA}

Any surjective equivalence is certainly an equivalence, so the `if'
direction is trivial. The force of the theorem is `only if'. 

Although Theorem~\ref{thm:A} is surely very old, it prepares us for similar
results in more complex contexts. The proof also sets the general
template. Given an equivalence of categories $F \from \A \to \B$, we wish
to construct a span $\A \ot \C \to \B$ of surjective equivalences, and we
do so by taking the universal span equipped with a natural isomorphism
\[
\begin{tikzcd}
\C \ar[d, "P" left] \ar[r, "Q"] 
\ar[rd, phantom,
"\rotatebox{45}{$\stackrel{\textstyle\Rightarrow}{\scriptstyle\sim}$}" description]        &
\B \ar[d, "1_\B"]  \\
\A \ar[r, "F" below]    &
\B.
\end{tikzcd}
\]
Concretely, an object of $\C$ consists of an object $a \in \A$, an object
$b \in \B$, and an isomorphism $F(a) \to b$. The maps in $\C$ are
as one would guess, so that $\C$ is a full subcategory of the comma
category $F \downarrow 1_\B$. The functors $P$ and $Q$ are the
projections, and both are surjective equivalences. 

Our second instance of Principle~\ref{princ:main} is as follows.

\begin{thmB}
Two monoidal categories $\A$ and $\B$ are monoidally equivalent if and only
if there exist a monoidal category $\C$ and strict monoidal surjective
equivalences $\A \ot \C \to \B$.
\end{thmB}

Here the word \emph{strict} refers to strict preservation of the monoidal
structure. For instance, if we write the span as $\A \otby{P} \C \toby{Q}
\B$ then the objects $P(c \otimes \twid{c})$ and $P(c) \otimes P(\twid{c})$
of $\A$ are \emph{equal} for all $c, \twid{c} \in \C$. The monoidal
categories involved are not strict, however.

The third instance concerns bicategories:

\begin{thmC}
Two bicategories $\A$ and $\B$ are biequivalent if and only if there exist
a bicategory $\C$ and strict surjective equivalences $\A \ot \C \to \B$. 
\end{thmC}

Again, \emph{strict} means strict preservation of composition in the
bicategories. A functor between bicategories is a \emph{surjective
equivalence} if it is genuinely surjective on objects and locally a
surjective equivalence of categories.

Generally, \emph{natural notion of equivalence} in
Principle~\ref{princ:main} means the canonical, weakest notion of
equivalence for the structures concerned, such as biequivalence for
bicategories. \emph{Strict} means strict preservation of the
structure. \emph{Surjective equivalence} means surjectivity on objects,
local surjectivity above dimension~$0$, and local injectivity in the top
dimension. It is a more primitive concept than ordinary equivalence, in
that it only refers to the underlying graphs of the categorical structures,
not the composition, monoidal product, etc.

(Incidentally, for $n$-dimensional structures, local injectivity in the top
dimension can be understood as local surjectivity in dimension $n +
1$. For example, a function between sets is injective if and only if the
corresponding functor between discrete categories is full.)

Principle~\ref{princ:main} states that the general notion of equivalence is
generated by the equivalences that are special in two respects: they are
\emph{literally} surjective in each dimension, and they \emph{strictly}
preserve the structure.

Some or all of Theorems~\ref{thm:A}, \ref{thm:B} and~\ref{thm:C} may be in
the literature, or folklore, but our final instance of
Principle~\ref{princ:main} appears to be new:

\begin{thmD}
Two double categories $\A$ and $\B$ are gregariously double equivalent if
and only if there exist a double category $\C$ and strict surjective
equivalences $\A \ot \C \to \B$.
\end{thmD}

This theorem suggests that whatever `gregariously double equivalent' means,
it is the \emph{right} notion of equivalence of double categories, since
the meaning of `strict surjective equivalence' is clear.  (The definitions
are stated in Section~\ref{sec:dbl}, and the aforementioned subtlety about
strictness is discussed there too.) Thus, Theorem~\ref{thm:D} supports the
conclusion of Campbell and of Moser, Sarazola and Verdugo.

All four proofs follow a similar pattern, suggested by the sketched proof
of Theorem~\ref{thm:A}. In all cases, given categorical structures $\A$ and
$\B$ and an equivalence $F \from \A \to \B$, we construct $\C$ by taking an
object to consist of an object $a \in \A$, an object $b \in \B$, and an
equivalence $F(a) \to b$ in $\B$. What `equivalence' means depends on the
context, and becomes more complex as we progress through the four
cases. For example, in the case of bicategories (Theorem~\ref{thm:C}), it
means not just a $1$-cell $F(a) \to b$ that is an equivalence in the usual
sense, but the full data of an adjoint equivalence between $F(a)$ and $b$.

\paragraph{Context} It remains to be seen how widely valid
Principle~\ref{princ:main} is. I learned of it from Carlos Simpson, who
used it in his work on Tamsamani's weak $n$-categories (or `$(n,
n)$-categories' in the current jargon). Our surjective equivalences are
what Simpson called easy equivalences. (See Section~2 and Corollary~6.6
of~\cite{SimpCMS}, as well as Pellissier~\cite{Pell} and Chapter~7 of
Simpson~\cite{SimpHTH}, which correct an error in~\cite{SimpCMS}.) Another
name for a surjective equivalence is a \emph{contractible} map, as in
Definition~2.3 of Cottrell~\cite{Cott} or Definition~P of~\cite{SDN}, for
instance.

Principle~\ref{princ:main} could be used to define equivalence between weak
$\infty$-categories (that is, $(\infty, \infty)$-categories) in situations
where no notion of weak $\infty$-functor is available. For example, three
different definitions of weak $\infty$-category proposed by Batanin, by
Leinster and by Penon (Definitions~B1, L1 and~P in~\cite{SDN}) each have
the following form: one constructs a certain monad $T$ on a certain
category of infinite-dimensional graphs, and one defines a weak
$\infty$-category to be a $T$-algebra, but the homomorphisms of
$T$-algebras should be understood as \emph{strict}
$\infty$-functors. Principle~\ref{princ:main} suggests declaring weak
$\infty$-categories $\A$ and $\B$ to be equivalent if and only if there
exists a span $\A \otby{P} \C \toby{Q} \B$ in the category of $T$-algebras
such that $P$ and $Q$ are locally surjective in every dimension.

The FOLDS system developed by Michael Makkai (`first order logic with
dependent sorts') also uses the idea that two structures should be
equivalent if and only if they are linked by a span of surjective
equivalences. This is his notion of $\mathbf{L}$-equivalence
(\cite{MakkFOL}, page~61). Probably Theorem~\ref{thm:A} also appears in
Makkai's FOLDS work. However, in the FOLDS setting, strictness of the
surjective equivalences is not a consideration.

The results in this paper, and the natural question `For which other
structures does Principle~\ref{princ:main} hold?', can be considered from
the perspective of model categories. Indeed, both the cited work of
Campbell and that of Moser, Sarazola and Verdugo are model category
theoretic, although they do not address Principle~\ref{princ:main}
directly. Work of Steve Lack puts a model structure on
bicategories~\cite{LackQMS}, and Principle~\ref{princ:main} is related
especially closely to the factorization lemma of Ken S.\ Brown
(\cite{BrowAHT}, page~421). But the approach taken here is wholly
elementary.

The term `surjective equivalence' was used for a slightly different
purpose by Blackwell, Kelly and Power (Remark~4.1 of~\cite{BKP}).

The structure of this paper is simple: Sections~\ref{sec:cats},
\ref{sec:mon}, \ref{sec:bi} and~\ref{sec:dbl} prove the theorems above for
categories, monoidal categories, bicategories and double categories,
respectively.

\section{Categories}
\label{sec:cats}

Here we state and prove the simple theorem on categories that is the
template for the later, more complex, results. Although the proof is
inevitably rather trivial, we spell it out nearly completely, since some
points that are straightforward for mere $1$-categories become substantial
in higher dimensions.

\begin{defn}
A functor between categories is a \demph{surjective equivalence} if it is
surjective on objects, full, and faithful.
\end{defn}

We will use the fact that categories $\A$ and $\B$ are equivalent if and
only if there exists a functor $\A \to \B$ that is essentially surjective
on objects, full, and faithful. In particular, a surjective equivalence is
an equivalence.

\begin{thm}
\label{thm:A}
Two categories $\A$ and $\B$ are equivalent if and only if there exist a
category $\C$ and surjective equivalences $\A \ot \C \to \B$.
\end{thm}

\begin{proof}
First suppose that such a span between $\A$ and $\B$ exists. Since a
surjective equivalence is an equivalence of categories, and equivalence of
categories is a symmetric transitive relation, $\A$ and $\B$ are
equivalent.

Conversely, suppose that $\A$ and $\B$ are equivalent. Take an equivalence
of categories $F \from \A \to \B$. Let $\C$ be the category whose objects
are triples $(a, b, \ell)$ consisting of an object $a \in \A$, an object $b
\in \B$, and an isomorphism $\ell \from Fa \to b$ in $\B$. A map $(a, b,
\ell) \to (a', b', \ell')$ in $\C$ is a pair $(f, g)$ consisting of a map
$f \from a \to a'$ in $\A$ and a map $g \from b \to b'$ in $\B$ such that
the square
\begin{equation}
\label{eq:cat-map}
\begin{tikzcd}
Fa \ar[r, "\ell", "\sim" below] \ar[d, "Ff" left]       &
b \ar[d, "g"]   \\
Fa' \ar[r, "\sim", "\ell'" below]       &
b'
\end{tikzcd}
\end{equation}
commutes. Composition and identities in $\C$ are defined in the obvious
way, so that the projections $\A \otby{P} \C \toby{Q} \B$ are functors.

It remains to show that $P$ and $Q$ are surjective equivalences. The proof
for $P$ does not need the hypothesis that $F$ is an equivalence, but the
proof for $Q$ does.

$P$ is surjective on objects: given $a \in \A$, the object $(a, Fa, 1_{Fa})$ of
$\C$ is mapped by $P$ to $a$.

$P$ is full: given $(a, b, \ell), (a', b', \ell') \in \C$ and $f \from a \to
a'$ in $\A$, put
\begin{equation}
\label{eq:cat-P-comp}
g = \Bigl( b \toby{\ell^{-1}} Fa \toby{Ff} Fa' \toby{\ell'} b' \Bigr).
\end{equation}
Then the square~\eqref{eq:cat-map} commutes, so $(f, g)$ is a map $(a, b,
\ell) \to (a', b', \ell')$ in $\C$ with $P(f, g) = f$. 

$P$ is faithful: in the fullness check just done, the commutativity
of~\eqref{eq:cat-map} forces $g$ to satisfy~\eqref{eq:cat-P-comp}.

This proves that $P \from \C \to \A$ is a surjective equivalence. Next we
show that $Q \from \C \to \B$ is too. 

$Q$ is surjective on objects: let $b \in \B$. Since $F$ is essentially
surjective on objects, we may choose an object $a \in \A$ and an
isomorphism $\ell \from Fa \to b$, and then $(a, b, \ell)$ is an object of
$\C$ with $Q(a, b, \ell) = b$.

$Q$ is full: given $(a, b, \ell), (a', b', \ell') \in \C$ and $g \from b
\to b'$ in $\B$, put
\begin{equation}
\label{eq:cat-Q-comp}
\hat{g} = 
\Bigl( Fa \toby{\ell} b \toby{g} b' \toby{\ell'^{-1}} Fa' \Bigr).
\end{equation}
Since $F$ is full, there is a map $f \from a \to a'$ in $\A$ such that $Ff
= \hat{g}$. Then the square~\eqref{eq:cat-map} commutes, so $(f, g)$ is a
map $(a, b, \ell) \to (a', b', \ell')$ in $\C$ with $Q(f, g) = g$. 

$Q$ is faithful: in the fullness check just done, any map $(\twid{f}, g)
\from (a, b, \ell) \to (a', b', \ell')$ in $\C$ must satisfy $F\twid{f} =
\hat{g}$ (by~\eqref{eq:cat-map}), and since $F$ is faithful, this
determines $\twid{f}$ uniquely.
\end{proof}

\section{Monoidal categories}
\label{sec:mon}

By definition, monoidal categories $\A$ and $\B$ (not necessarily strict)
are \demph{monoidally equivalent} if there exist strong monoidal functors
$F \from \A \oppairu \B \from G$ and monoidal natural isomorphisms $1_\A
\toby{\sim} GF$ and $FG \toby{\sim} 1_\B$. A logically equivalent condition
is that there exists a strong monoidal functor $\A \to \B$ that is an
equivalence of categories.

\begin{thm}
\label{thm:B}
Two monoidal categories $\A$ and $\B$ are monoidally equivalent if and only
if there exist a monoidal category $\C$ and strict monoidal surjective
equivalences $\A \ot \C \to \B$.
\end{thm}

\begin{proof}
`If' follows from the facts just stated about monoidal equivalence.

For `only if', take a monoidal equivalence $(F, \phi) \from \A \to \B$,
where the symbol $\phi$ stands for the family of coherence isomorphisms
\[
\bigl(
\phi_{a, \twid{a}} \from 
Fa \otimes F\twid{a} 
\toby{\sim} 
F(a \otimes \twid{a}) 
\bigr)_{a, \twid{a} \in \A}
\]
together with the coherence isomorphism $\phi_\cdot \from I_\B \toby{\sim}
FI_\A$.  
 
We build on the proof of Theorem~\ref{thm:A}. Define the category $\C$ and
the surjective equivalences $\A \otby{P} \C \toby{Q} \B$ as in that
proof. Then $\C$ can be given a monoidal structure as follows. On objects,
the tensor product is defined by
\[
(a, b, \ell) \otimes (\twid{a}, \twid{b}, \twid{\ell}) 
=
\Bigl( 
a \otimes \twid{a}, \ 
b \otimes \twid{b}, \ 
F(a \otimes \twid{a}) \toby{\phi_{a, \twid{a}}^{-1}}
Fa \otimes F\twid{a} \toby{\ell \otimes \twid{\ell}} 
b \otimes \twid{b}
\Bigr),
\]
and on maps, 
\[
(f, g) \otimes (\twid{f}, \twid{g}) 
= 
(f \otimes \twid{f}, g \otimes \twid{g}).
\]
The unit object is $(I_\A, I_\B, FI_\A \toby{\phi_\cdot^{-1}} I_\B)$. The
associativity and unit isomorphisms of $\C$ are inherited from those of
$\A$ and $\B$, as are the coherence axioms.

By construction, the surjective equivalences $P$ and $Q$ in the proof of
Theorem~\ref{thm:A} are strict monoidal functors, and the result follows.
\end{proof}

\section{Bicategories}
\label{sec:bi}

To prove the theorem on bicategories, we first recall some facts
about both equivalence \emph{in} a bicategory and equivalence
\emph{between} bicategories. 

Objects $b$ and $b'$ of a bicategory $\B$ are \demph{equivalent}, written
as $b \eqv b'$, if there exist 1-cells $\oppair{b}{b'}{\ell}{r}$ such that
$r \ell \iso 1_b$ and $\ell r \iso 1_{b'}$. An equivalent condition is that
there exists an \demph{adjoint equivalence} between $b$ and $b'$, that is,
not only $1$-cells $\ell$ and $r$ as above, but also invertible $2$-cells
\[
\eta \from 1_b \toby{\sim} r\ell,
\quad
\epsilon \from \ell r \toby{\sim} 1_{b'}
\]
satisfying the triangle identities for an adjunction.

Now let $\A$ and $\B$ be bicategories. Among the various classes of maps
between bicategories, we will work with the pseudofunctors (which preserve
composition of $1$-cells up to coherent isomorphism) and the strict
$2$-functors (which preserve composition strictly). 

The pseudofunctors $\A \to \B$, together with the pseudonatural
isomorphisms and modifications between them, form the functor bicategory
$[\A, \B]$. By definition, $\A$ and $\B$ are \demph{biequivalent} if and
only if there exist pseudofunctors $\oppair{\A}{\B}{F}{G}$ such that $GF
\eqv 1_\A$ in $[\A, \A]$ and $FG \eqv 1_\B$ in $[\B, \B]$. A logically
equivalent condition is that there exists a pseudofunctor $F \from \A \to
\B$ that is \demph{essentially surjective on objects} (for all $b \in B$,
there exists $a \in \A$ such that $Fa \eqv b$) and locally an equivalence
of categories. (See, for example, Section~9 of Street~\cite{StreCS}.)

\begin{defn}
A pseudofunctor $F \from \A \to \B$ is a \demph{surjective equivalence} if it
is surjective on objects (for all $b \in \B$, there exists $a \in \A$ such
that $Fa = b$) and locally a surjective equivalence of categories.
\end{defn}

Thus, a surjective equivalence of bicategories is a pseudofunctor that is
surjective on objects, locally surjective on objects, and locally full and
faithful. This is a stronger condition than being a biequivalence.

\begin{thm}
\label{thm:C}
Two bicategories $\A$ and $\B$ are biequivalent if and only if there exist
a bicategory $\C$ and strict surjective equivalences $\A \ot \C \to \B$.
\end{thm}

This result does not subsume Theorem~\ref{thm:C}, since even if $\A$ and
$\B$ have only one object, $\C$ in general has more. Nor is it subsumed by
Theorem~\ref{thm:D} on double categories, for similar reasons.

For simplicity, the proof will be written as if $\A$ and $\B$ were strict
2-categories, as is justified by the coherence theorem. Routine checks are
omitted.

\begin{proof}
For `if', biequivalence is an equivalence relation, so if such surjective
equivalences exist then $\A$, $\C$ and $\B$ are biequivalent.

For `only if', let $F \from \A \to \B$ be a biequivalence. Define a
bicategory $\C$ as follows.

An object $(a, b, \ell, r, \eta, \epsilon)$ of $\C$ consists of an object
$a \in \A$, an object $b \in \B$, and $1$- and $2$-cells
\[
\oppair{Fa}{b,}{\ell}{r}
\quad
1_{Fa} \toby{\eta} r\ell,
\quad
\ell r \toby{\epsilon} 1_b
\]
defining an adjoint equivalence between $Fa$ and $b$.

A $1$-cell 
\[
(f, g, \lambda, \rho) \from 
(a, b, \ell, r, \eta, \epsilon) \to (a', b', \ell', r', \eta', \epsilon')
\]
in $\C$ consists of $1$-cells
\[
a \toby{f} a', 
\quad
b \toby{g} b'
\]
in $\A$ and $\B$ respectively, together with invertible $2$-cells
\[
\begin{tikzcd}[row sep = tiny]
&
Fa' \ar[rd, "\ell'"] 
\ar[dd, phantom, "\sim\Downarrow\lambda" description]  &       
\\
Fa \ar[ru, "Ff"] \ar[rd, "\ell"']       &
&
b'      \\
&
b \ar[ru, "g"']
\end{tikzcd}
\qquad
\begin{tikzcd}[row sep = tiny]
&
Fa \ar[rd, "Ff"] 
\ar[dd, phantom, "\sim\Downarrow\rho" description]  &       
\\
b \ar[ru, "r"] \ar[rd, "g"']       &
&
Fa'     \\
&
b' \ar[ru, "r'"']
\end{tikzcd}
\]
in $\B$ satisfying coherence axioms that we state using string diagrams:
\begin{equation}
\label{eq:bi-coh-1}
\begin{doublestring}
\node[ghost] (eta'1) {};
\node[ghost] (eta'2) [right=of eta'1] {};
\node[ghost] (eta'out) [left=of eta'1] {};
\node[w2] (eta') at ($(eta'1)!1/2!(eta'2)$) {$\eta'$};
\node[ghost] (lambda2) [below=of eta'1] {};
\node[ghost] (lambda1) [left=of lambda2] {};
\node[ghost] (lambdaout) [right=of lambda2] {};
\node[w2] (lambda) at ($(lambda1)!1/2!(lambda2)$) {$\lambda$};
\inedge{eta'1}{lambda2}{\ell'}
\outabovex{lambda1}{eta'out}{Ff}
\outbelow{lambda1}{\ell}
\outbelow{lambda2}{g}
\outbelowx{eta'2}{lambdaout}{r'}
\end{doublestring}
=
\begin{doublestring}
\node[ghost] (eta1) {};
\node[ghost] (eta2) [right=of eta1] {};
\node[ghost] (etaout) [right=of eta2] {};
\node[w2] (eta) at ($(eta1)!1/2!(eta2)$) {$\eta$};
\node[ghost] (rhoout) [below=of eta1] {};
\node[ghost] (rho1) [right=of rhoout] {};
\node[ghost] (rho2) [right=of rho1] {};
\node[w2] (rho) at ($(rho1)!1/2!(rho2)$) {$\rho$};
\inedge{eta2}{rho1}{r}
\outabovex{rho2}{etaout}{Ff}
\outbelowx{eta1}{rhoout}{\ell}
\outbelow{rho1}{g}
\outbelow{rho2}{r'}
\end{doublestring}
\quad\quad
\begin{doublestring}
\node[ghost] (lambda1) {};
\node[ghost] (lambda2) [right=of lambda1] {};
\node[ghost] (lambdaout) [left=of lambda1] {};
\node[w2] (lambda) at ($(lambda1)!1/2!(lambda2)$) {$\lambda$};
\node[ghost] (epsilon2) [below=of lambda1] {};
\node[ghost] (epsilon1) [left=of epsilon2] {};
\node[ghost] (epsilonout) [right=of epsilon2] {};
\node[w2] (epsilon) at ($(epsilon1)!1/2!(epsilon2)$) {$\epsilon$};
\inedge{lambda1}{epsilon2}{\ell}
\outabovex{epsilon1}{lambdaout}{r}
\outabove{lambda1}{Ff}
\outabove{lambda2}{\ell'}
\outbelowx{lambda2}{epsilonout}{g}
\end{doublestring}
=
\begin{doublestring}
\node[ghost] (rho1) {};
\node[ghost] (rho2) [right=of rho1] {};
\node[ghost] (rhoout) [right=of rho2] {};
\node[w2] (rho) at ($(rho1)!1/2!(rho2)$) {$\rho$};
\node[ghost] (epsilon'out) [below=of eta1] {};
\node[ghost] (epsilon'1) [right=of epsilon'out] {};
\node[ghost] (epsilon'2) [right=of epsilon'1] {};
\node[w2] (epsilon') at ($(epsilon'1)!1/2!(epsilon'2)$) {$\epsilon'$};
\inedge{rho2}{epsilon'1}{r'}
\outabove{rho1}{r}
\outabove{rho2}{Ff}
\outabovex{epsilon'2}{rhoout}{\ell'}
\outbelowx{rho1}{epsilon'out}{g}
\end{doublestring}
\end{equation}

A $2$-cell
\[
\begin{tikzcd}[column sep = huge]
(a, b, \ell, r, \eta, \epsilon)         
\ar[r, bend left, "{(f, g, \lambda, \rho)}"]
\ar[r, bend right, "{(\twid{f}, \twid{g}, \twid{\lambda}, \twid{\rho})}"']
\ar[r, phantom, "{\Downarrow(\alpha, \beta)}" description]       
&
(a', b', \ell', r', \eta', \epsilon')
\end{tikzcd}
\]
in $\C$ consists of $2$-cells
\[
\begin{tikzcd}[column sep = large]
a 
\ar[r, bend left, "f"] \ar[r, bend right, "\twid{f}"']
\ar[r, phantom, "\Downarrow\alpha" description] &
a'
\end{tikzcd}
\qquad
\begin{tikzcd}[column sep = large]
b
\ar[r, bend left, "g"] \ar[r, bend right, "\twid{g}"']
\ar[r, phantom, "\Downarrow\beta" description] &
b'
\end{tikzcd}
\]
in $\A$ and $\B$ respectively, satisfying the coherence axioms
\begin{equation}
\label{eq:bi-coh-2}
\begin{doublestring}
\node[basic] (Falpha) {$F\alpha$};
\node[ghost] (Falphaout) [right=of Falpha] {};
\node[ghost] (lambdatwid1) [below=of Falpha] {};
\node[ghost] (lambdatwid2) [right=of lambdatwid1] {};
\node[w2] (lambdatwid) at ($(lambdatwid1)!1/2!(lambdatwid2)$) {\raisebox{-1mm}{$\twid{\lambda}$}};
\inedge{Falpha}{lambdatwid1}{F\twid{f}}
\outabove{Falpha}{Ff}
\outabovex{lambdatwid2}{Falphaout}{\ell'}
\outbelow{lambdatwid1}{\ell}
\outbelow{lambdatwid2}{\twid{g}}
\end{doublestring}
=
\begin{doublestring}
\node[ghost] (lambda1) {};
\node[ghost] (lambda2) [right=of lambda1] {};
\node[w2] (lambda) at ($(lambda1)!1/2!(lambda2)$) {$\lambda$};
\node[basic] (beta) [below=of lambda2] {$\beta$};
\node[ghost] (betaout) [left=of beta] {};
\inedges{lambda2}{beta}{g}
\outabove{lambda1}{Ff}
\outabove{lambda2}{\ell'}
\outbelowx{lambda1}{betaout}{\ell}
\outbelow{beta}{\twid{g}}
\end{doublestring}
\qquad
\begin{doublestring}
\node[basic] (Falpha) {\makebox[0mm]{$F\alpha$}};
\node[ghost] (Falphaout) [left=of Falpha] {};
\node[ghost] (rhotwid1) [below=of Falphaout] {};
% \node[ghost] (rhotwid2) [right=of rhotwid1] {};
\node[ghost] (rhotwid2) [right=of rhotwid1] {};
\node[w2] (rhotwid) at ($(rhotwid1)!1/2!(rhotwid2)$) {$\twid{\rho}$};
\inedges{Falpha}{rhotwid2}{F\twid{f}}
\outabove{Falpha}{Ff}
\outabovex{rhotwid1}{Falphaout}{r}
\outbelow{rhotwid1}{\twid{g}}
\outbelow{rhotwid2}{r'}
\end{doublestring}
=
\begin{doublestring}
\node[ghost] (rho1) {};
\node[ghost] (rho2) [right=of rho1] {};
\node[w2] (rho) at ($(rho1)!1/2!(rho2)$) {$\rho$};
\node[basic] (beta) [below=of rho1] {$\beta$};
\node[ghost] (betaout) [right=of beta] {};
\inedge{rho1}{beta}{g}
\outabove{rho1}{r}
\outabove{rho2}{Ff}
\outbelow{beta}{\twid{g}}
\outbelowx{rho2}{betaout}{r'}
\end{doublestring}
\end{equation}
Composition, identities and coherence cells in $\C$ are defined in the
obvious way, using the coherence cells of the pseudofunctor $F$ as
necessary. The evident projections
\[
\A \otby{P} \C \toby{Q} \B
\]
are then strict $2$-functors.

We briefly digress on some redundancy in the definition of $\C$. In the
definition of $1$-cell, the coherence axioms~\eqref{eq:bi-coh-1} state that
$\lambda$ and $\rho$ are mates under the adjunctions $\ell \ladj r$ and
$\ell' \ladj r'$ (Section~2 of Kelly and Street~\cite{KeStRE2}). Thus,
$\lambda$ and $\rho$ determine each other; for example,
\begin{equation}
\label{eq:bi-mate}
\rho
=
\begin{doublestring}
\node[ghost] (eta'1) {};
\node[ghost] (eta'2) [right=of eta'1] {};
\node[ghost] (eta'out2) [left=of eta'1] {};
\node[ghost] (eta'out1) [left=of eta'out2] {};
\node[w2] at ($(eta'1)!1/2!(eta'2)$) {$\eta'$};
\node[ghost] (lambda1) [below=of eta'out2] {};
\node[ghost] (lambda2) [right=of lambda1] {};
\node[w2] at ($(lambda1)!1/2!(lambda2)$) {$\lambda$};
\node[ghost] (ep2) [below=of lambda1] {};
\node[ghost] (ep1) [left=of ep2] {};
\node[ghost] (epout1) [right=of ep2] {};
\node[ghost] (epout2) [right=of epout1] {};
\node[w2] at ($(ep1)!1/2!(ep2)$) {$\epsilon$};
\inedge{eta'1}{lambda2}{\ell'}
\inedges{lambda1}{ep2}{\ell}
\outabovex{ep1}{eta'out1}{r}
\outabovex{lambda1}{eta'out2}{Ff}
\outbelowx{lambda2}{epout1}{g}
\outbelowx{eta'2}{epout2}{r'}
\end{doublestring}
\end{equation}
So a $1$-cell $(a, b, \ldots) \to (a', b', \ldots)$ in $\C$ could
equivalently be defined as consisting of $1$-cells $f \from a \to a'$ in
$\A$ and $g \from b \to b'$ in $\B$ together with an invertible $2$-cell
$\lambda \from \ell' \of Ff \to g \of \ell$ in $\B$ whose mate is also
invertible. The mateship just noted also means that in the definition of
$2$-cell, each of the equations~\eqref{eq:bi-coh-2} implies the other. For
example, the first implies the second by~\eqref{eq:bi-mate} and its
analogue for $\twid{\lambda}$ and $\twid{\rho}$. But it seems preferable to
maintain the symmetry.

Resuming the proof, it remains to show that both $P$ and $Q$ are surjective
equivalences. As in Theorems~\ref{thm:A} and~\ref{thm:B}, the hypothesis
that $F$ is a biequivalence is used in the proof for $Q$ but not $P$. We
begin with $P$.

$P$ is surjective on objects: given $a \in \A$, the object $Fa$ of $\B$
together with the identity adjoint equivalence on $Fa$ defines an object of
$\C$ whose image under $P$ is equal to $a$.

$P$ is locally surjective on objects: let $(a, b, \ell, r, \eta, \epsilon),
(a', b', \ldots) \in \C$ and $f \from a \to a'$ in $\A$. Define
\[
g = \Bigl( 
b \toby{r} Fa \toby{Ff} Fa' \toby{\ell'} b'
\Bigr)
\]
(a $1$-cell in $\B$), then put 
\[
\lambda
=
\begin{doublestring}
\node[ghost] (eta1) {};
\node[ghost] (eta2) [right=of eta1] {};
\node[ghost] (etaout1) [right=of eta2] {};
\node[ghost] (etaout2) [right=of etaout1] {};
\node[w2] at ($(eta1)!1/2!(eta2)$) {$\eta$};
\node[ghost] (e1) [below=of eta2] {};
\node[ghost] (e2) [right=of e1] {};
\node[ghost] (e3) [right=of e2] {};
\node[ghost] (eout) [left=of e1] {};
\node[w3] at ($(e1)!1/2!(e3)$) {$=$};
\inedge{eta2}{e1}{r}
\outabovex{e2}{etaout1}{Ff}
\outabovex{e3}{etaout2}{\ell'}
\outbelowx{eta1}{eout}{\ell}
\outbelow{e2}{g}
\end{doublestring}
\qquad
\rho
=
\begin{doublestring}
\node[ghost] (eta'1) {};
\node[ghost] (eta'2) [right=of eta'1] {};
\node[ghost] (eta'out2) [left=of eta'1] {};
\node[ghost] (eta'out1) [left=of eta'out2] {};
\node[w2] at ($(eta'1)!1/2!(eta'2)$) {$\eta'$};
\node[ghost] (e3) [below=of eta'1] {};
\node[ghost] (e2) [left=of e3] {};
\node[ghost] (e1) [left=of e2] {};
\node[ghost] (eout) [right=of e3] {};
\node[w3] at (e2) {$=$};
\inedge{eta'1}{e3}{\ell'}
\outabovex{e1}{eta'out1}{r}
\outabovex{e2}{eta'out2}{Ff}
\outbelow{e2}{g}
\outbelowx{eta'2}{eout}{r'}
\end{doublestring}
\]
(both invertible $2$-cells in $\B$). The coherence
axioms~\eqref{eq:bi-coh-1} in the definition of $1$-cell in $\C$ follow, in
this case, from the triangle identities for $(\eta, \epsilon)$ and $(\eta',
\epsilon')$. Hence $(f, g, \lambda, \rho)$ is a $1$-cell $(a, b, \ldots)
\to (a', b', \ldots)$ in $\C$, with $P(f, g, \lambda, \rho) = f$.

$P$ is locally full and faithful: let 
\begin{equation}
\label{eq:bi-empty}
\begin{tikzcd}[column sep = huge]
(a, b, \ell, r, \eta, \epsilon)         
\ar[r, bend left, "{(f, g, \lambda, \rho)}"]
\ar[r, bend right, "{(\twid{f}, \twid{g}, \twid{\lambda}, \twid{\rho})}"']
&
(a', b', \ell', r', \eta', \epsilon')
\end{tikzcd}
\end{equation}
in $\C$ and let 
\[
\begin{tikzcd}[column sep = large]
a 
\ar[r, bend left, "f"] \ar[r, bend right, "\twid{f}"']
\ar[r, phantom, "\Downarrow\alpha" description] &
a'
\end{tikzcd}
\]
in $\A$. We must show that there is a unique $2$-cell
\[
\begin{tikzcd}[column sep = large]
b
\ar[r, bend left, "g"] \ar[r, bend right, "\twid{g}"']
\ar[r, phantom, "\Downarrow\beta" description] &
b'
\end{tikzcd}
\]
in $\B$ such that $(\alpha, \beta)$ is a $2$-cell $(f, g, \lambda, \rho)
\to (\twid{f}, \twid{g}, \twid{\lambda}, \twid{\rho})$ in $\C$, that is,
such that equations~\eqref{eq:bi-coh-2} hold. We have already seen that the
second of equations~\eqref{eq:bi-coh-2} is redundant, and since $\lambda$
is invertible, the first is equivalent to 
\begin{equation}
\label{eq:bi-ff}
\begin{doublestring}
\node[basic] (beta) {$\beta$};
\node[ghost] (betaout) [left=of beta] {};
\outabove{beta}{g}
\outbelow{beta}{\twid{g}}
\throughvert{betaout}{betaout}{\ell}
\end{doublestring}
=
\begin{doublestring}
\node[ghost] (linv1) {};
\node[ghost] (linv2) [right=of linv1] {};
\node[w2] at ($(linv1)!1/2!(linv2)$) {$\lambda^{-1}$};
\node[basic] (Falpha) [below=of linv1] {$F\alpha$};
\node[ghost] (ltwid1) [below=of Falpha] {};
\node[ghost] (ltwid2) [right=of ltwid1] {};
\node[w2] at ($(ltwid1)!1/2!(ltwid2)$) {\raisebox{-1mm}{$\twid{\lambda}$}};
\inedges{linv1}{Falpha}{} %{Ff}
\inedges{Falpha}{ltwid1}{} %{F\twid{f}}
\inedge{linv2}{ltwid2}{} %{\ell'}
\outabove{linv1}{\ell}
\outabove{linv2}{g}
\outbelow{ltwid1}{\ell}
\outbelow{ltwid2}{\twid{g}}
\end{doublestring}
\end{equation}
We have to show that~\eqref{eq:bi-ff} has a unique solution $\beta$. This
can be done by an elementary argument, but a little bicategory theory makes
it easier. The representable $\B(-, b') \from \B^\op \to \Cat$, like all
pseudofunctors, preserves equivalences. Since $\ell$ is an equivalence, so
is $\B(\ell, b')$; that is, composition with $\ell$ defines an equivalence
of categories
\[
- \of \ell \from \B(b, b') \to \B(Fa, b').
\]
In particular, $- \of \ell$ is full and faithful, so~\eqref{eq:bi-ff}
has a unique solution $\beta$. 

This completes the proof that the first projection $P \from \C \to \A$ is a
surjective equivalence. We now show that the second projection $Q \from \C
\to \B$ is too. 

$Q$ is surjective on objects: let $b \in \B$. Since $F$ is essentially
surjective on objects, we can choose an object $a \in \A$ and an adjoint
equivalence $(\ell, r, \eta, \epsilon)$ between $Fa$ and $b$. Then $(a, b,
\ell, r, \eta, \epsilon)$ is an object of $\C$ whose image under $Q$ is
equal to $b$.

$Q$ is locally surjective on objects: let $(a, b, \ell, r, \eta, \epsilon),
(a', b', \ldots) \in \C$ and $g \from b \to b'$ in $\B$. Since $F$ is
locally essentially surjective on objects, there exist a $1$-cell $f \from a
\to a'$ in $\A$ and an invertible $2$-cell
\[
\begin{tikzcd}[row sep = 0.3ex]
Fa \ar[rrr, "Ff", bend left = 15]    
\ar[rrr, phantom, "\sim\Downarrow\chi" description]
\ar[rd, "\ell"']
&
&
&
Fa'     \\
&
b \ar[r, "g"']  &
b' \ar[ru, "r'"']
\end{tikzcd}
\]
in $\B$. Define $2$-cells $\lambda$ and $\rho$ in $\B$ by
\[
\lambda
=
\begin{doublestring}
\node[ghost] (chi1) {};
\node[ghost] (chi2) [right=of chi1] {};
\node[ghost] (chi3) [right=of chi2] {};
\node[ghost] (chiout) [right=of chi3] {};
\node[w3] at (chi2) {$\chi$};
\node[ghost] (ep'1) [below=of chi3] {};
\node[ghost] (ep'2) [right=of ep'1] {};
\node[ghost] (ep'out2) [left=of ep'1] {};
\node[ghost] (ep'out1) [left=of ep'out2] {};
\node[w2] at ($(ep'1)!1/2!(ep'2)$) {$\epsilon'$};
\inedge{chi3}{ep'1}{r'}
\outabove{chi2}{Ff}
\outabovex{ep'2}{chiout}{\ell'}
\outbelowx{chi1}{ep'out1}{\ell}
\outbelowx{chi2}{ep'out2}{g}
\end{doublestring}
\qquad
\rho
=
\begin{doublestring}
\node[ghost] (chi1) {};
\node[ghost] (chi2) [right=of chi1] {};
\node[ghost] (chi3) [right=of chi2] {};
\node[ghost] (chiout) [left=of chi1] {};
\node[w3] at (chi2) {$\chi$};
\node[ghost] (ep2) [below=of chi1] {};
\node[ghost] (ep1) [left=of ep2] {};
\node[ghost] (epout1) [right=of ep2] {};
\node[ghost] (epout2) [right=of epout1] {};
\node[w2] at ($(ep1)!1/2!(ep2)$) {$\epsilon$};
\inedge{chi1}{ep2}{\ell}
\outabovex{ep1}{chiout}{r}
\outabove{chi2}{Ff}
\outbelowx{chi2}{epout1}{g}
\outbelowx{chi3}{epout2}{r'}
\end{doublestring}
\]
Both $\lambda$ and $\rho$ are invertible, since $\chi$, $\epsilon$ and
$\epsilon'$ are. The coherence axioms~\eqref{eq:bi-coh-1} are easily
checked. Thus, $(f, g, \lambda, \rho)$ is a $1$-cell $(a, b, \ldots) \to
(a', b', \ldots)$ in $\C$, with $Q(f, g, \lambda, \rho) = g$.

$Q$ is locally full and faithful: the proof of this is very similar to the
proof that $P$ is locally full and faithful, now also using the hypothesis
that $F$ is locally full and faithful. 
\end{proof}

\section{Double categories}
\label{sec:dbl}

For simplicity, we restrict to \emph{strict} double categories. We use
string diagrams: for example, 2-cells
\[
\begin{tikzcd}
a \ar[r, "f"] \ar[d, "s"']       
\ar[dr, phantom, "\alpha" description]      
&
b \ar[d, "t"]   \\
c \ar[r, "g"']  &
d
\end{tikzcd}
\qquad
\begin{tikzcd}
a \ar[rr, equal] \ar[d, equal]
\ar[rrd, phantom, "\eta" description]    &
&
a \ar[d, equal] \\
a \ar[r, "\ell"']       &
b \ar[r, "r"']  &
a
\end{tikzcd}
\]
in a double category, where the elongated equality signs denote identity
cells, are drawn as 
\[
\begin{doublestring}
\node[basic] (alpha) {$\alpha$};
\outabove{alpha}{f}
\outright{alpha}{t}
\outbelow{alpha}{g}
\outleft{alpha}{s}
\end{doublestring}
\qquad
\begin{doublestring}
\node[ghost] (eta1) {};
\node[ghost] (eta2) [right=of eta1]{};
\node[w2] at ($(eta1)!1/2!(eta2)$) {$\eta$};
\outbelow{eta1}{\ell}
\outbelow{eta2}{r}
\end{doublestring}
\]
(These are similar to the string diagrams used by Arkor and
McDermott~\cite{ArMc}, but different from those of Myers~\cite{MyerSDD}.)

A particular family of examples of double categories will illuminate some
of the definitions that follow.

\begin{example}
\label{eg:2-dbl}
Given a 2-category $\scat{K}$, there is a double category
$\scat{D}(\scat{K})$ whose objects are those of $\scat{K}$, whose
horizontal and vertical $1$-cells are both the $1$-cells of $\scat{K}$, and
whose $2$-cells
\[
\begin{tikzcd}
a \ar[r, "f"] \ar[d, "s"']       
\ar[dr, phantom, "\alpha" description]      
&
b \ar[d, "t"]   \\
c \ar[r, "g"']  &
d
\end{tikzcd}
\]
are the $2$-cells
\[
\begin{tikzcd}
a \ar[r, "f"] \ar[d, "s"']       
\ar[dr, phantom, "\swnt\alpha" description]      
&
b \ar[d, "t"]   \\
c \ar[r, "g"']  &
d
\end{tikzcd}
\]
in $\scat{K}$ (that is, $\alpha \from t \of f \to g \of s$). Composition is
defined in the obvious way.
\end{example}

We now recall a concept due to Grandis and Par\'e (\cite{GrPaADC},
Section~1.2).

\begin{defn}
\label{defn:companion}
Two $1$-cells $a \toby{f} b$ and $\vertmap{a}{b}{f'}$ in a double category
are \demph{companions} if there exist $2$-cells
\[
\begin{tikzcd}
a \ar[r, equal] \ar[d, equal] 
\ar[rd, phantom, "\sigma" description]     &
a \ar[d, "f'"]  \\
a \ar[r, "f"']  &
b
\end{tikzcd}
\qquad
\begin{tikzcd}
a \ar[r, "f"] \ar[d, "f'"'] 
\ar[rd, phantom, "\tau" description] &
b \ar[d, equal] \\
b \ar[r, equal] &
b
\end{tikzcd}
\]
satisfying the equations
\begin{equation}
\label{eq:companion}
\begin{doublestring}
\node[basic] (tau) {$\tau$};
\node[basic] (sigma) [left=of tau] {$\sigma$};
\inedge{sigma}{tau}{f'}
\outabove{tau}{f}
\outbelow{sigma}{f}
\end{doublestring}
=
\begin{doublestring}
\node[ghost] (anchor) {};
\throughvert{anchor}{anchor}{f}
\end{doublestring}
\qquad
\begin{doublestring}
\node[basic] (tau) {$\tau$};
\node[basic] (sigma) [above=of tau] {$\sigma$};
\inedge{sigma}{tau}{f}
\outright{sigma}{f'}
\outleft{tau}{f'}
\end{doublestring}
=
\begin{doublestring}
\node[ghost] (anchor) {};
\throughhoriz{anchor}{anchor} {f'}
\end{doublestring}
\end{equation}
Following Dawson, Par\'e and Pronk, we call $\sigma$ and $\tau$
\demph{binding cells} (Definition~3.4 of~\cite{DPP}).
\end{defn}

\begin{example}
In the double category $\scat{D}(\scat{K})$ derived from a $2$-category
$\scat{K}$, the $1$-cells $f$ and $f'$ are companions if and only if they
are isomorphic in $\scat{K}$. In an arbitrary double category, it makes no
sense to ask whether a horizontal $1$-cell and a vertical $1$-cell are
isomorphic, but this example suggests that companionship is a kind of
substitute.
\end{example}

We now build up to Campbell's definition of gregarious double equivalence
between double categories~\cite{Camp}, beginning with his notion of
equivalence between objects \emph{in} a double category.

\begin{defn}
Objects $a$ and $b$ of a double category are \demph{gregariously
equivalent} if there exist companion $1$-cells $a \toby{f} b$ and
$\vertmap{a}{b}{f'}$ such that $f$ is a horizontal equivalence and $f'$ is
a vertical equivalence.
\end{defn}

Thus, for $a$ and $b$ to be gregariously equivalent, it is not enough for
them to be horizontally equivalent and vertically equivalent in unrelated
ways. There must be horizontal and vertical equivalences that are related
by companionship. It is a pleasant exercise to check that gregarious
equivalence is an equivalence relation.

(The first use of the word \emph{gregarious} in a related context appears
to have been by Dawson, Par\'e and Pronk; see Definitions~2.4 and~3.11
of~\cite{DPP}.)

To discuss equivalence \emph{between} double categories, we first need the
right notion of functor. A \demph{double pseudofunctor} is a map between
double categories that preserves both horizontal and vertical $1$-cell
composition up to coherent isomorphism, and $2$-cell composition
strictly. Thus, a double pseudofunctor $F \from \A \to \B$ consists of
not only the expected assignations on $0$-, $1$- and $2$-cells, but also
vertically invertible coherence cells
\[
\begin{tikzcd}
Fa \ar[r, "Ff"] \ar[d, equal]
\ar[rrd, phantom, "\phi_{f, f'}" description]   &
Fa' \ar[r, "Ff'"]       &
Fa'' \ar[d, equal]      \\
Fa \ar[rr, "F(f' \of f)"']      &
&
Fa''
\end{tikzcd}
\qquad
\begin{tikzcd}
Fa \ar[r, equal] \ar[d, equal]
\ar[rd, phantom, "\phi_a" description]  &
Fa \ar[d, equal]        \\
Fa \ar[r, "F(1_a)"']    &
Fa
\end{tikzcd}
\]
(for each composable pair $a \toby{f} a' \toby{f'} a''$ of horizontal
$1$-cells in $\A$ in the first case, and each object $a \in \A$ in the
second), and, similarly, horizontally invertible coherence $2$-cells for
vertical composition of $1$-cells, all subject to coherence axioms. Details
are given in Definition~6.1 of Shulman~\cite{ShulCCL}. A double
pseudofunctor whose coherence cells are all identities is a \demph{strict
double functor}.

\begin{defn}
Let $\A$ and $\B$ be double categories. A double pseudofunctor $F \from \A
\to \B$ is a \demph{gregarious double equivalence} (and $\A$ and $\B$ are
\demph{gregariously double equivalent}) if:
\begin{itemize}
\item 
$F$ is \demph{surjective on objects up to gregarious equivalence}: for
every $b \in \B$, there exists $a \in \A$ such that $Fa$ is gregariously
equivalent to $b$;

\item
$F$ is \demph{horizontally essentially full}: for each $a, a' \in \A$ and
horizontal $1$-cell $Fa \toby{g} Fa'$ in $\B$, there exist a horizontal
$1$-cell $a \toby{f} a'$ in $\A$ and a vertically invertible $2$-cell
\[
\begin{tikzcd}
Fa \ar[r, "Ff"] \ar[d, equal] 
\ar[rd, phantom, "\chi" description]     &
Fa' \ar[d, equal]       \\
Fa \ar[r, "g"']        &
Fa';
\end{tikzcd}
\]

\item
$F$ is \demph{vertically essentially full}, similarly;

\item
$F$ is \demph{full and faithful on $2$-cells}: for every diagram
\[
\begin{tikzcd}
a \ar[r, "f"] \ar[d, "s"']      &a' \ar[d, "s'"]        \\
\twid{a} \ar[r, "\twid{f}"']    &\twid{a}'
\end{tikzcd}
\]
of $0$- and $1$-cells in $\A$ and every $2$-cell
\[
\begin{tikzcd}
Fa \ar[r, "Ff"] \ar[d, "Fs"']      
\ar[rd, phantom, "\beta" description]   &
Fa' \ar[d, "Fs'"]        \\
F\twid{a} \ar[r, "F\twid{f}"']     &F\twid{a}'
\end{tikzcd}
\]
in $\B$, there exists a unique $2$-cell
\[
\begin{tikzcd}
a \ar[r, "f"] \ar[d, "s"']      
\ar[rd, phantom, "\alpha" description]  &
a' \ar[d, "s'"]        \\
\twid{a} \ar[r, "\twid{f}"']    &\twid{a}'
\end{tikzcd}
\]
in $\A$ such that $F\alpha = \beta$.
\end{itemize}
\end{defn}

It is straightforward to check that gregarious double equivalence is a
reflexive, transitive relation on the class of double categories.  Symmetry
is less straightforward, but will follow as a corollary of our main
theorem.

\begin{remark}
The double functors in Campbell's slides~\cite{Camp} are implicitly strict,
but for our purposes, it is essential to use the pseudo version. Otherwise,
the relation of gregarious double equivalence is not symmetric.

To see the problem, first consider monoidal categories $\A$ and $\B$. It
may be that there exists a strict monoidal equivalence $\A \to \B$ but no
strict monoidal equivalence $\B \to \A$, even if $\A$ and $\B$ themselves
are strict. For example, let $\B$ be the discrete category on the
$2$-element group $\Z/2\Z$, with monoidal structure given by addition. Let
$\A$ be the category whose object-set is $\Z$, with one map $m \to n$ when
$m \equiv n \pmod{2}$ and none otherwise, again made into a monoidal
category by addition. Let $F \from \A \to \B$ be the usual quotient
map. Then $F$ is a strict monoidal functor between strict monoidal
categories, and it is a surjective equivalence. However, there is no strict
monoidal equivalence $G \from \B \to \A$. For since $1 + 1 = 0$ in
$\Z/2\Z$, if we put $n = G(1) \in \Z$ then $n + n = 0$ in $\Z$, so $n = 0$;
hence $G(1) = G(0)$, contradicting $G$ being an equivalence. 

By viewing a strict monoidal category as a degenerate double category, we
obtain an example of double categories $\A$ and $\B$ such that there exists
a strict double functor $\A \to \B$ that is a gregarious double
equivalence, but no such functor $\B \to \A$. 
\end{remark}

\begin{defn}
A double pseudofunctor $F \from \A \to \B$ is a \demph{surjective
equivalence} if:
\begin{itemize}
\item 
$F$ is surjective on objects;

\item
$F$ is \demph{horizontally full}: for each $a, a' \in \A$ and horizontal
$1$-cell $Fa \toby{g} Fa'$ in $\B$, there exists a horizontal $1$-cell $a
\toby{f} a'$ in $\A$ such that $Ff = g$;

\item
$F$ is \demph{vertically full}, similarly;

\item
$F$ is full and faithful on $2$-cells.
\end{itemize}
\end{defn}

It is immediate that any surjective equivalence is a gregarious double
equivalence. Our main theorem is:

\begin{thm}
\label{thm:D}
Two double categories $\A$ and $\B$ are gregariously double equivalent if
and only if there exist a double category $\C$ and strict surjective
equivalences $\A \ot \C \to \B$.
\end{thm}

Here, a \demph{strict surjective equivalence} is a strict double functor
that is a surjective equivalence.

We now prepare for the proof of Theorem~\ref{thm:D}.

\begin{lemma}
\label{lemma:inverse}
Let $\C$ and $\A$ be double categories. If there exists a strict surjective
equivalence $\C \to \A$ then there exists a gregarious double equivalence
$\A \to \C$.
\end{lemma}

The lemma is true under the weaker hypothesis that there exists a gregarious
double equivalence $\C \to \A$ (Corollary~\ref{cor:sym}), but this is all we
will need.

\begin{proof}
Let $P \from \C \to \A$ be a strict surjective equivalence. We will
define a gregarious double equivalence $J \from \A \to \C$.

On objects: for each $a \in \A$, choose $Ja \in \C$ such that $PJa = a$,
which is possible since $P$ is surjective on objects.

On horizontal $1$-cells: for each $a \toby{f} a'$ in $\A$, choose a
horizontal $1$-cell $Ja \toby{Jf} Ja'$ in $\B$ such that $PJf = f$, which
is possible since $P$ is horizontally full.

On vertical $1$-cells, define $J$ similarly.

On $2$-cells: for each $2$-cell
\[
\begin{tikzcd}
a \ar[r, "f"] \ar[d, "s"'] 
\ar[rd, phantom, "\alpha" description]  &
a' \ar[d, "s'"]         \\
\twid{a} \ar[r, "\twid{f}"']    &
\twid{a}'
\end{tikzcd}
\]
in $\A$, there is a unique $2$-cell
\[
\begin{tikzcd}
Ja \ar[r, "Jf"] \ar[d, "Js"'] 
\ar[rd, phantom, "J\alpha" description]  &
Ja' \ar[d, "Js'"]         \\
J\twid{a} \ar[r, "J\twid{f}"']    &
J\twid{a}'
\end{tikzcd}
\]
in $\C$ such that $PJ\alpha = \alpha$, since $P$ is full and faithful on
$2$-cells. 

To define the coherence isomorphisms of the double pseudofunctor $J$, first
take a composable pair $a \toby{f} a' \toby{f'} a''$ of horizontal
$1$-cells in $\A$. Then $PJ(f' \of f) = f' \of f$ by definition of $J$, but
also 
\[
P(J(f') \of J(f)) = PJ(f') \of PJ(f) = f' \of f
\]
because $P$ is strict and by definition of $J$ again. Since $P$ is full and
faithful on $2$-cells, there is a unique $2$-cell
\[
\begin{tikzcd}
Ja \ar[r, "Jf"] \ar[d, equal]
\ar[rrd, phantom, "\zeta_{f, f'}" description] &
Ja' \ar[r, "Jf'"]       &
Ja'' \ar[d, equal]      \\
Ja \ar[rr, "J(f') \of J(f)"']   &
&
Ja''
\end{tikzcd}
\]
in $\C$ such that $P(\zeta_{f, f'}) = 1_{f' \of f}$. The other coherence
cells are defined similarly, and it is readily checked that they are
invertible and satisfy the coherence axioms for a double pseudofunctor.

This defines the double pseudofunctor $J \from \A \to \C$. By construction,
$P \of J$ is equal to the identity on $\A$. It remains to show that $J$ is
a gregarious double equivalence.

$J$ is surjective on objects up to gregarious double equivalence: let $c
\in \C$. It is enough prove that $JPc$ and $c$ are gregariously
equivalent. 

Since $P(JPc) = (PJ)Pc = Pc$ and $P$ is horizontally full, we can choose a
horizontal $1$-cell $JPc \toby{h} c$ such that $Ph = 1_{Pc}$. Choose a
vertical $1$-cell $\vertmap{JPc}{c}{h'}$ similarly. Now, the image under
$P$ of the unfilled square
\[
\begin{tikzcd}
JPc \ar[r, equal] \ar[d, equal]   &
JPc \ar[d, "h'"]        \\
JPc \ar[r, "h"']        &
c
\end{tikzcd}
\]
is
\[
\begin{tikzcd}
Pc \ar[r, equal] \ar[d, equal]  &Pc \ar[d, equal]       \\
Pc \ar[r, equal]                &Pc,
\end{tikzcd}
\]
and since $P$ is full and faithful on $2$-cells, there is a unique $2$-cell 
\[
\begin{tikzcd}
JPc \ar[r, equal] \ar[d, equal]   
\ar[rd, phantom, "\sigma" description]  &
JPc \ar[d, "h'"]        \\
JPc \ar[r, "h"']        &
c
\end{tikzcd}
\]
such that $P\sigma$ is the identity $2$-cell on $Pc$. Similarly, there is a
unique $2$-cell
\[
\begin{tikzcd}
JPc \ar[r, "h"] \ar[d, "h'"']   
\ar[rd, phantom, "\tau" description]    &
c \ar[d, equal]        \\
c \ar[r, equal]                 &c
\end{tikzcd}
\]
such that $P\tau$ is the identity. With $\sigma$ and $\tau$ as binding
cells, $h$ and $h'$ are companions; indeed, the fact that $\sigma$ and
$\tau$ satisfy equations~\eqref{eq:companion} follows from $P$ being full
and faithful and $P\sigma$ and $P\tau$ being identities. (More generally, a
double functor that is full and faithful on $2$-cells always reflects
companionship.)

Moreover, since the strict double functor $P \from \C \to \A$ is a
surjective equivalence, the induced $2$-functor on the underlying
horizontal $2$-categories is also a surjective equivalence, and in
particular, a biequivalence. Since $Ph$ is an equivalence in the horizontal
$2$-category of $\A$, it follows that $h$ is an equivalence in the
horizontal $2$-category of $\C$. Similarly, $h'$ is a vertical
equivalence. Hence $JPc$ and $c$ are gregariously equivalent, completing
the proof that $J$ is a surjective on objects up to gregarious equivalence.

$J$ is horizontally essentially full: let $a, a' \in \A$ and $Ja \toby{h}
Ja'$ in $\C$. Since $P \of J = 1_{\A}$, we also have the horizontal
$1$-cell $Ja \toby{JPh} Ja'$ in $\C$, and we will prove that $JPh \iso h$
in the underlying horizontal $2$-category of $\C$.

Indeed, $PJPh = Ph$, and $P$ is full and faithful on $2$-cells, so there is
a unique $2$-cell 
\[
\begin{tikzcd}
Ja \ar[r, "JPh"] \ar[d, equal]
\ar[rd, phantom, "\chi" description]  &
Ja' \ar[d, equal]       \\
Ja \ar[r, "h"']         &
Ja'
\end{tikzcd}
\]
in $\C$ such that $P\chi = 1_{Ph}$. Again because $P$ is full and
faithful on $2$-cells, $\chi$ is vertically invertible. Hence $JPh \iso
h$, as required.

$J$ is vertically essentially full, similarly.

$J$ is full and faithful on $2$-cells: this follows easily from the same
property of $P$ together with the fact that $PJ = 1_\A$.

This completes the proof that $J$ is a gregarious double equivalence.
\end{proof}

The proof of Theorem~\ref{thm:D} will use the notions of horizontal and
vertical adjunction in a double category. A \demph{horizontal
adjunction} in a double category consists of horizontal $1$-cells
$\oppair{a}{b}{\ell}{r}$ and $2$-cells
\[
\begin{doublestring}
\node[ghost] (eta1) {};
\node[ghost] (eta2) [right=of eta1] {};
\node[w2] at ($(eta1)!1/2!(eta2)$) {$\eta$};
\outbelow{eta1}{\ell}
\outbelow{eta2}{r}
\end{doublestring}
\qquad
\begin{doublestring}
\node[ghost] (ep1) {};
\node[ghost] (ep2) [right=of ep1] {};
\node[w2] at ($(ep1)!1/2!(ep2)$) {$\epsilon$};
\outabove{ep1}{r}
\outabove{ep2}{\ell}
\end{doublestring}
\]
satisfying the familiar triangle identities. A \demph{vertical
adjunction} consists of vertical $1$-cells $\begin{tikzcd}[cramped,
sep=small] a \ar[d, shift right, "\ell"']\\ b \ar[u, shift right, "r"']
\end{tikzcd}$ and $2$-cells
\[
\begin{doublestring}
\node[ghost] (eta1) {};
\node[ghost] (eta2) [below=of eta1] {};
\node[h2] at ($(eta1)!1/2!(eta2)$) {$\eta$};
\outleft{eta1}{\ell}
\outleft{eta2}{r}
\end{doublestring}
\qquad
\begin{doublestring}
\node[ghost] (ep1) {};
\node[ghost] (ep2) [below=of ep1] {};
\node[h2] at ($(ep1)!1/2!(ep2)$) {$\epsilon$};
\outright{ep1}{r}
\outright{ep2}{\ell}
\end{doublestring}
\]
(note the orientations!)\ satisfying triangle identities. The orientations
are chosen so that in the double category $\scat{D}(\scat{K})$ of
Example~\ref{eg:2-dbl}, an adjunction of either of these two types amounts
to an adjunction $\ell \ladj r$ in the original $2$-category $\scat{K}$. 

A \demph{horizontal adjoint equivalence} is a horizontal adjunction in
which $\eta$ and $\epsilon$ are vertically invertible, and similarly
\demph{vertical adjoint equivalences}.

Although the proof of Theorem~\ref{thm:D} has a certain notational
complexity, the main idea is simple and already apparent in the proof of
the analogous Theorem~\ref{thm:A} on categories. There, given an
equivalence of categories $F \from \A \to \B$, a new category was
constructed not by taking objects to be pairs $(a \in \A, b \in \B)$
\emph{such that} $Fa \iso b$, but by taking them to be pairs $(a, b)$
together with a \emph{specified} isomorphism $Fa \to b$. 

Similarly, for bicategories (Theorem~\ref{thm:C}), an object of $\C$
involved not just a choice of $a$ and $b$, and not just a choice of
$1$-cell $Fa \toby{\ell} b$ for which a pseudoinverse \emph{exists}, but a
\emph{specified} pseudoinverse $b \toby{r} Fa$ and \emph{specified}
isomorphisms $\eta\from 1 \to r\ell$ and $\epsilon \from \ell r \to 1$,
together with the imposition of all reasonable coherence axioms (in this
case, the triangle identities). 

The same principle applies to double categories. The length of the proof is
due only to the necessity of handling more data, not because of greater
conceptual depth. Like many proofs in two-dimensional category theory, much
of the work is clerical and has a strong sense of inevitability about it,
making one dream of an approach in which the drudgework is unnecessary.

\begin{pfof}{Theorem~\ref{thm:D}}
Let $\A$ and $\B$ be double categories. First suppose that there exist
strict surjective equivalences $\A \otby{P} \C \toby{Q} \B$. By
Lemma~\ref{lemma:inverse}, there exists a gregarious double equivalence $J
\from \A \to \C$. Since $Q$ is also a gregarious double equivalence, so is
the composite $QJ \from \A \to \B$. Hence $\A$ and $\B$ are gregariously
double equivalent.

Conversely, suppose there exists a gregarious double equivalence $F \from
\A \to \B$. Define a double category $\C$ as follows.

An object of $\C$ consists of an object $a \in \A$, an object $b \in \B$, a
horizontal adjoint equivalence
\[
\begin{tikzcd}
Fa \ar[r, "\ell_\Ho", shift left] &
b \ar[l, "r_\Ho", shift left]
\end{tikzcd}
\qquad
\begin{doublestring}
\node[ghost] (eta1) {};
\node[ghost] (eta2) [right=of eta1] {};
\node[w2] at ($(eta1)!1/2!(eta2)$) {$\eta_\Ho$};
\outbelow{eta1}{\ell_\Ho}
\outbelow{eta2}{r_\Ho}
\end{doublestring}
\qquad
\begin{doublestring}
\node[ghost] (ep1) {};
\node[ghost] (ep2) [right=of ep1] {};
\node[w2] at ($(ep1)!1/2!(ep2)$) {$\epsilon_\Ho$};
\outabove{ep1}{r_\Ho}
\outabove{ep2}{\ell_\Ho}
\end{doublestring}
\]
in $\B$, a vertical adjoint equivalence
\[
\begin{tikzcd}
Fa \ar[d, shift right, "\ell_\Ve"'] \\
b \ar[u, shift right, "r_\Ve"']
\end{tikzcd}
\qquad
\begin{doublestring}
\node[ghost] (eta1) {};
\node[ghost] (eta2) [below=of eta1] {};
\node[h2] at ($(eta1)!1/2!(eta2)$) {$\eta_\Ve$};
\outleft{eta1}{\ell_\Ve}
\outleft{eta2}{r_\Ve}
\end{doublestring}
\qquad
\begin{doublestring}
\node[ghost] (ep1) {};
\node[ghost] (ep2) [below=of ep1] {};
\node[h2] at ($(ep1)!1/2!(ep2)$) {$\epsilon_\Ve$};
\outright{ep1}{r_\Ve}
\outright{ep2}{\ell_\Ve}
\end{doublestring}
\]
in $\B$, and $2$-cells
\[
\begin{doublestring}
\node[basic] (sigma) {$\sigma$};
\outright{sigma}{\ell_\Ve}
\outbelow{sigma}{\ell_\Ho}
\end{doublestring}
\qquad
\begin{doublestring}
\node[basic] (tau) {$\tau$};
\outabove{tau}{\ell_\Ho}
\outleft{tau}{\ell_\Ve}
\end{doublestring}
\]
making $\ell_\Ho$ and $\ell_\Ve$ into companions (as in
Definition~\ref{defn:companion}). We write such an object as $(a, b, \ell,
r, \eta, \epsilon, \sigma, \tau)$.

(This definition may appear asymmetric, in that it mentions companionship
between the left adjoints but not the right. We could make it look more
symmetric by including $2$-cells $\ovln{\sigma}$ and $\ovln{\tau}$ making
$r_\Ho$ and $r_\Ve$ into companions, then imposing appropriate coherence
axioms relating $\ovln{\sigma}$ and $\ovln{\tau}$ to $\sigma$ and
$\tau$. In fact, this approach is equivalent to the one taken here. For
instance, $\ovln{\sigma}$ and $\ovln{\tau}$ are determined by $\sigma$ and
$\tau$ as follows:
\begin{equation}
\label{eq:dbl-mates}
\ovln{\sigma}
=
\begin{doublestring}
\node[ghost] (eta1) {};
\node[ghost] (eta2) [below=of eta1] {};
\node[ghost] (eta3) [below=of eta2] {};
\node[h3] at (eta2) {$\eta_\Ve$};
\node[basic] (sigma) [left=of eta1] {$\sigma$};
\node[ghost] (sigmaout) [left=of sigma] {};
\node[ghost] (ep2) [below=of sigma] {};
\node[ghost] (ep1) [left=of ep2] {};
\node[ghost] (epout) [below=of ep1] {};
\node[w2] at ($(ep1)!1/2!(ep2)$) {$\epsilon_\Ho$};
\inedges{sigma}{eta1}{\ell_\Ve}
\inedges{sigma}{ep2}{\ell_\Ho\!}
\outabovex{ep1}{sigmaout}{r_\Ho}
\outleftx{eta3}{epout}{r_\Ve}
\end{doublestring}
\qquad
\ovln{\tau}
=
\begin{doublestring}
\node[ghost] (ep1) {};
\node[ghost] (ep2) [below=of ep1] {};
\node[ghost] (ep3) [below=of ep2] {};
\node[h3] at (ep2) {$\epsilon_\Ve$};
\node[basic] (tau) [right=of ep3] {$\tau$};
\node[ghost] (tauout) [right=of tau] {};
\node[ghost] (eta1) [above=of tau] {};
\node[ghost] (eta2) [right=of eta1] {};
\node[ghost] (etaout) [above=of eta2] {};
\node[w2] at ($(eta1)!1/2!(eta2)$) {$\eta_\Ho$};
\inedge{ep3}{tau}{\ell_\Ve}
\inedge{eta1}{tau}{\ell_\Ho}
\outrightx{ep1}{etaout}{r_\Ve}
\outbelowx{eta2}{tauout}{r_\Ho}
\end{doublestring}
\end{equation}
Here we sacrifice symmetry to reduce the notational load.)

A horizontal $1$-cell 
\[
(f, g, \lambda, \rho) \from
(a, b, \ell, r, \eta, \epsilon, \sigma, \tau) \to
(a', b', \ell', r', \eta', \epsilon', \sigma', \tau') 
\]
in $\C$ consists of horizontal $1$-cells
\[
a \toby{f} a',
\qquad
b \toby{g} b'
\]
in $\A$ and $\B$ respectively, invertible $2$-cells
\[
\begin{tikzcd}
Fa \ar[r, "Ff"] \ar[d, equal]
\ar[rrd, phantom, "\lambda_\Ho" description]    &
Fa' \ar[r, "\ell'_\Ho"] &
b' \ar[d, equal]        \\
Fa \ar[r, "\ell_\Ho"']  &
b \ar[r, "g"']  &
b'
\end{tikzcd}
\qquad
\begin{tikzcd}
b \ar[r, "r_\Ho"] \ar[d, equal]
\ar[rrd, phantom, "\rho_\Ho" description]       &
Fa \ar[r, "Ff"] &
Fa' \ar[d, equal]       \\
b \ar[r, "g"']  &
b' \ar[r, "r'_\Ho"']    &
Fa'
\end{tikzcd}
\]
in $\B$, and $2$-cells
\[
\begin{tikzcd}
Fa \ar[r, "Ff"] \ar[d, "\ell_\Ve"'] 
\ar[rd, phantom, "\lambda_\Ve" description]     &
Fa' \ar[d, "\ell'_\Ve"] \\
b \ar[r, "g"']  &
b'
\end{tikzcd}
\qquad
\begin{tikzcd}
b \ar[r, "g"] \ar[d, "r_\Ve"']
\ar[rd, phantom, "\rho_\Ve" description]        &
b' \ar[d, "r'_\Ve"]      \\
Fa \ar[r, "Ff"'] &
Fa'
\end{tikzcd}
\]
in $\B$, satisfying the following coherence axioms. First, $\lambda_\Ho$
and $\rho_\Ho$ are compatible (much as in the proof of Theorem~\ref{thm:C},
equation~\eqref{eq:bi-coh-1}):
\begin{equation}
\label{eq:dbl-coh-1ha}
\begin{doublestring}
\node[ghost] (eta'1) {};
\node[ghost] (eta'2) [right=of eta'1] {};
\node[ghost] (eta'out) [left=of eta'1] {};
\node[w2] (eta') at ($(eta'1)!1/2!(eta'2)$) {$\eta'_\Ho$};
\node[ghost] (lambda2) [below=of eta'1] {};
\node[ghost] (lambda1) [left=of lambda2] {};
\node[ghost] (lambdaout) [right=of lambda2] {};
\node[w2] (lambda) at ($(lambda1)!1/2!(lambda2)$) {$\lambda_\Ho$};
\inedge{eta'1}{lambda2}{\ell'_\Ho}
\outabovex{lambda1}{eta'out}{Ff}
\outbelow{lambda1}{\ell_\Ho}
\outbelow{lambda2}{g}
\outbelowx{eta'2}{lambdaout}{r'_\Ho}
\end{doublestring}
=
\begin{doublestring}
\node[ghost] (eta1) {};
\node[ghost] (eta2) [right=of eta1] {};
\node[ghost] (etaout) [right=of eta2] {};
\node[w2] (eta) at ($(eta1)!1/2!(eta2)$) {$\eta_\Ho$};
\node[ghost] (rhoout) [below=of eta1] {};
\node[ghost] (rho1) [right=of rhoout] {};
\node[ghost] (rho2) [right=of rho1] {};
\node[w2] (rho) at ($(rho1)!1/2!(rho2)$) {$\rho_\Ho$};
\inedge{eta2}{rho1}{r_\Ho}
\outabovex{rho2}{etaout}{Ff}
\outbelowx{eta1}{rhoout}{\ell_\Ho}
\outbelow{rho1}{g}
\outbelow{rho2}{r'_\Ho}
\end{doublestring}
\quad\quad
\begin{doublestring}
\node[ghost] (lambda1) {};
\node[ghost] (lambda2) [right=of lambda1] {};
\node[ghost] (lambdaout) [left=of lambda1] {};
\node[w2] (lambda) at ($(lambda1)!1/2!(lambda2)$) {$\lambda_\Ho$};
\node[ghost] (epsilon2) [below=of lambda1] {};
\node[ghost] (epsilon1) [left=of epsilon2] {};
\node[ghost] (epsilonout) [right=of epsilon2] {};
\node[w2] (epsilon) at ($(epsilon1)!1/2!(epsilon2)$) {$\epsilon_\Ho$};
\inedge{lambda1}{epsilon2}{\ell_\Ho}
\outabovex{epsilon1}{lambdaout}{r_\Ho}
\outabove{lambda1}{Ff}
\outabove{lambda2}{\ell'_\Ho}
\outbelowx{lambda2}{epsilonout}{g}
\end{doublestring}
=
\begin{doublestring}
\node[ghost] (rho1) {};
\node[ghost] (rho2) [right=of rho1] {};
\node[ghost] (rhoout) [right=of rho2] {};
\node[w2] (rho) at ($(rho1)!1/2!(rho2)$) {$\rho_\Ho$};
\node[ghost] (epsilon'out) [below=of eta1] {};
\node[ghost] (epsilon'1) [right=of epsilon'out] {};
\node[ghost] (epsilon'2) [right=of epsilon'1] {};
\node[w2] (epsilon') at ($(epsilon'1)!1/2!(epsilon'2)$) {$\epsilon'_\Ho$};
\inedge{rho2}{epsilon'1}{r'_\Ho}
\outabove{rho1}{r_\Ho}
\outabove{rho2}{Ff}
\outabovex{epsilon'2}{rhoout}{\ell'_\Ho}
\outbelowx{rho1}{epsilon'out}{g}
\end{doublestring}
\end{equation}
Second, $\lambda_\Ve$ and $\rho_\Ve$ are compatible:
\begin{equation}
\label{eq:dbl-coh-1hb}
\begin{doublestring}
\node[basic] (lambda) {$\lambda_\Ve$};
\node[basic] (rho) [below=of lambda] {$\rho_\Ve$};
\node[ghost] (eta'1) [right=of lambda] {};
\node[ghost] (eta'2) [right=of rho] {};
\node[h2] at ($(eta'1)!1/2!(eta'2)$) {$\eta'_\Ve$};
\inedges{lambda}{rho}{g}
\inedge{lambda}{eta'1}{\ell'_\Ve}
\inedges{rho}{eta'2}{r'_\Ve}
\outabove{lambda}{Ff}
\outleft{lambda}{\ell_\Ve}
\outleft{rho}{r_\Ve}
\outbelow{rho}{Ff}
\end{doublestring}
=
\begin{doublestring}
\node[ghost] (eta1) {};
\node[ghost] (eta2) [below=of eta1] {};
\node[ghost] (etaout1) [right=of eta1] {};
\node[ghost] (etaout2) [right=of eta2] {};
\node[h2] at ($(eta1)!1/2!(eta2)$) {$\eta_\Ve$};
\throughvert{etaout1}{etaout2}{Ff}
\outleft{eta1}{\ell_\Ve}
\outleft{eta2}{r_\Ve}
\end{doublestring}
\qquad
\begin{doublestring}
\node[basic] (rho) {$\rho_\Ve$};
\node[basic] (lambda) [below=of rho] {$\lambda_\Ve$};
\node[ghost] (ep1) [left=of rho] {};
\node[ghost] (ep2) [left=of lambda] {};
\node[h2] at ($(ep1)!1/2!(ep2)$) {$\epsilon_\Ve$};
\inedge{rho}{lambda}{Ff}
\inedge{ep1}{rho}{r_\Ve}
\inedges{ep2}{lambda}{\ell_\Ve}
\outabove{rho}{g}
\outright{rho}{r'_\Ve}
\outright{lambda}{\ell'_\Ve}
\outbelow{lambda}{g}
\end{doublestring}
=
\begin{doublestring}
\node[ghost] (ep'1) {};
\node[ghost] (ep'2) [below=of ep'1] {};
\node[ghost] (ep'out1) [left=of ep'1] {};
\node[ghost] (ep'out2) [left=of ep'2] {};
\node[h2] at ($(ep'1)!1/2!(ep'2)$) {$\epsilon'_\Ve$};
\throughvert{ep'out1}{ep'out2}{g}
\outright{ep'1}{r'_\Ve}
\outright{ep'2}{\ell'_\Ve}
\end{doublestring}
\end{equation}
Third and finally, $\lambda_\Ho$ and $\lambda_\Ve$ are compatible:
\begin{equation}
\label{eq:dbl-coh-1hc}
\begin{doublestring}
\node[ghost] (lambda1) {};
\node[ghost] (lambda2) [right=of lambda1] {};
\node[w2] at ($(lambda1)!1/2!(lambda2)$) {$\lambda_\Ho$};
\node[basic] (sigma') [above=of lambda2] {$\sigma'$};
\node[ghost] (sigma'out) [left=of sigma'] {};
\inedges{sigma'}{lambda2}{\ell'_\Ho}
\outabovex{lambda1}{sigma'out}{Ff}
\outright{sigma'}{\ell'_\Ve}
\outbelow{lambda2}{g}
\outbelow{lambda1}{\ell_\Ho}
\end{doublestring}
=
\begin{doublestring}
\node[basic] (lambda) {$\lambda_\Ve$};
\node[basic] (sigma) [left=of lambda] {$\sigma$};
\inedge{sigma}{lambda}{\ell_\Ve}
\outabove{lambda}{Ff}
\outright{lambda}{\ell'_\Ve}
\outbelow{lambda}{g}
\outbelow{sigma}{\ell_\Ho}
\end{doublestring}
\qquad
\begin{doublestring}
\node[ghost] (lambda1) {};
\node[ghost] (lambda2) [right=of lambda1] {};
\node[w2] at ($(lambda1)!1/2!(lambda2)$) {$\lambda_\Ho$};
\node[basic] (tau) [below=of lambda1] {$\tau$};
\node[ghost] (tauout) [right=of tau] {};
\inedge{lambda1}{tau}{\ell_\Ho}
\outabove{lambda1}{Ff}
\outabove{lambda2}{\ell'_\Ho}
\outbelowx{lambda2}{tauout}{g}
\outleft{tau}{\ell_\Ve}
\end{doublestring}
=
\begin{doublestring}
\node[basic] (lambda) {$\lambda_\Ve$};
\node[basic] (tau') [right=of lambda] {$\tau'$};
\inedge{lambda}{tau'}{\ell'_\Ve}
\outabove{lambda}{Ff}
\outabove{tau'}{\ell'_\Ho}
\outbelow{lambda}{g}
\outleft{lambda}{\ell_\Ve}
\end{doublestring}
\end{equation}
(The analogous compatibility condition for $\rho_\Ho$ and $\rho_\Ve$,
involving the cells $\ovln{\sigma}$ and $\ovln{\tau}$ of
equation~\eqref{eq:dbl-mates}, follows from
conditions~\eqref{eq:dbl-coh-1ha}--\eqref{eq:dbl-coh-1hc}.) 

Vertical $1$-cells in $\C$ are, of course, defined similarly, but we spell
it out. A vertical $1$-cell
\[
\begin{tikzcd}
(a, b, \ell, r, \eta, \epsilon, \sigma, \tau)     
\ar[d, "{(s, t, \gamma, \delta)}"]        \\
(\twid{a}, \twid{b}, \twid{\ell}, \twid{r}, 
\twid{\eta}, \twid{\epsilon}, \twid{\sigma}, \twid{\tau})     
\end{tikzcd}
\]
in $\C$ consists of vertical $1$-cells
\[
\begin{tikzcd}
a \ar[d, "s"] \\
\twid{a}
\end{tikzcd}
\qquad
\begin{tikzcd}
b \ar[d, "t"] \\
\twid{b}
\end{tikzcd}
\]
in $\A$ and $\B$ respectively, horizontally invertible $2$-cells 
\[
\begin{tikzcd}
Fa \ar[r, equal] \ar[d, "\ell_\Ve"']
\ar[rdd, phantom, "\gamma_\Ve" description]     &
Fa \ar[d, "Fs"] \\
b \ar[d, "t"']  &
F\twid{a} \ar[d, "\twid{\ell}_\Ve"]     \\
\twid{b}' \ar[r, equal]         &
\twid{b}'
\end{tikzcd}
\qquad
\begin{tikzcd}
b \ar[r, equal] \ar[d, "t"']
\ar[rdd, phantom, "\delta_\Ve" description]     &
b \ar[d, "r_\Ve"]       \\
\twid{b} \ar[d, "\twid{r}_\Ve"']        &
Fa \ar[d, "Fs"] \\
F\twid{a} \ar[r, equal]         &
F\twid{a}
\end{tikzcd}
\]
in $\B$, and $2$-cells
\[
\begin{tikzcd}
Fa \ar[r, "\ell_\Ho"] \ar[d, "Fs"']     
\ar[rd, phantom, "\gamma_\Ho" description]      &
b \ar[d, "t"]   \\
F\twid{a} \ar[r, "\twid{\ell}_\Ho"']    &
\twid{b}
\end{tikzcd}
\qquad
\begin{tikzcd}
b \ar[r, "r_\Ho"] \ar[d, "t"']
\ar[rd, phantom, "\delta_\Ho" description]      &
Fa \ar[d, "Fs"] \\
\twid{b} \ar[r, "\twid{r}_\Ho"']        &
F\twid{a}
\end{tikzcd}
\]
in $\B$, satisfying the following coherence axioms. First, $\gamma_\Ve$ and
$\delta_\Ve$ are compatible:
\begin{equation}
\label{eq:dbl-coh-1va}
\begin{doublestring}
\node[ghost] (gamma1) {};
\node[ghost] (gamma2) [below=of gamma1] {};
\node[ghost] (gammaout) [below=of gamma2] {};
\node[h2] at ($(gamma1)!1/2!(gamma2)$) {$\gamma_\Ve$};
\node[ghost] (etat1) [right=of gamma2] {};
\node[ghost] (etat2) [below=of etat1] {};
\node[ghost] (etatout) [above=of etat1] {};
\node[h2] at ($(etat1)!1/2!(etat2)$) {$\twid{\eta}_\Ve$};
\inedge{gamma2}{etat1}{\twid{\ell}_\Ve}
\outrightx{gamma1}{etatout}{Fs}
\outleftx{etat2}{gammaout}{\twid{r}_\Ve}
\outleft{gamma2}{t}
\outleft{gamma1}{\ell_\Ve}
\end{doublestring}
=
\begin{doublestring}
\node[ghost] (eta1) {};
\node[ghost] (eta2) [below=of eta1] {};
\node[ghost] (etaout) [below=of eta2] {};
\node[h2] at ($(eta1)!1/2!(eta2)$) {$\eta_\Ve$};
\node[ghost] (delta1) [left=of gamma2] {};
\node[ghost] (delta2) [below=of delta1] {};
\node[ghost] (deltaout) [above=of delta1] {};
\node[h2] at ($(delta1)!1/2!(delta2)$) {$\delta_\Ve$};
\inedge{delta1}{eta2}{r_\Ve}
\outrightx{delta2}{etaout}{Fs}
\outleft{delta2}{\twid{r}_\Ve}
\outleft{delta1}{t}
\outleftx{eta1}{deltaout}{\ell_\Ve}
\end{doublestring}
\qquad
\begin{doublestring}
\node[ghost] (epsilon1) {};
\node[ghost] (epsilon2) [below=of epsilon1] {};
\node[ghost] (epsilonout) [below=of epsilon2] {};
\node[h2] at ($(epsilon1)!1/2!(epsilon2)$) {$\epsilon_\Ve$};
\node[ghost] (gamma1) [right=of epsilon2] {};
\node[ghost] (gamma2) [below=of gamma1] {};
\node[ghost] (gammaout) [above=of gamma1] {};
\node[h2] at ($(gamma1)!1/2!(gamma2)$) {$\gamma_\Ve$};
\inedge{epsilon2}{gamma1}{\ell_\Ve}
\outrightx{epsilon1}{gammaout}{r_\Ve}
\outright{gamma1}{Fs}
\outright{gamma2}{\twid{\ell}_\Ve}
\outleftx{gamma2}{epsilonout}{t}
\end{doublestring}
=
\begin{doublestring}
\node[ghost] (delta1) {};
\node[ghost] (delta2) [below=of delta1] {};
\node[ghost] (deltaout) [below=of delta2] {};
\node[h2] at ($(delta1)!1/2!(delta2)$) {$\delta_\Ve$};
\node[ghost] (ept1) [left=of delta2] {};
\node[ghost] (ept2) [below=of ept1] {};
\node[ghost] (eptout) [above=of ept1] {};
\node[h2] at ($(ept1)!1/2!(ept2)$) {$\twid{\epsilon}_\Ve$};
\inedge{ept1}{delta2}{\twid{r}_\Ve}
\outright{delta1}{r_\Ve}
\outright{delta2}{Fs}
\outrightx{ept2}{deltaout}{\twid{\ell}_\Ve}
\outleftx{delta1}{eptout}{t}
\end{doublestring}
\end{equation}
Second, $\gamma_\Ho$ and $\delta_\Ho$ are compatible:
\begin{equation}
\label{eq:dbl-coh-1vb}
\begin{doublestring}
\node[basic] (gamma) {$\gamma_\Ho$};
\node[basic] (delta) [right=of gamma] {$\delta_\Ho$};
\node[ghost] (eta1) [above=of gamma] {};
\node[ghost] (eta2) [above=of delta] {};
\node[w2] at ($(eta1)!1/2!(eta2)$) {$\eta_\Ho$};
\inedges{eta1}{gamma}{\ell_\Ho}
\inedge{eta2}{delta}{r_\Ho}
\inedges{gamma}{delta}{t}
\outright{delta}{Fs}
\outbelow{delta}{\twid{r}_\Ho}
\outbelow{gamma}{\twid{\ell}_\Ho}
\outleft{gamma}{Fs}
\end{doublestring}
=
\begin{doublestring}
\node[ghost] (etat1) {};
\node[ghost] (etat2) [right=of etat1] {};
\node[ghost] (etatout1) [above=of etat1] {};
\node[ghost] (etatout2) [above=of etat2] {};
\node[w2] at ($(etat1)!1/2!(etat2)$) {$\twid{\eta}_\Ho$};
\throughhoriz{etatout1}{etatout2}{Fs}
\outbelow{etat1}{\twid{\ell}_\Ho}
\outbelow{etat2}{\twid{r}_\Ho}
\end{doublestring}
\qquad
\begin{doublestring}
\node[basic] (delta) {$\delta_\Ho$};
\node[basic] (gamma) [right=of delta] {$\gamma_\Ho$};
\node[ghost] (epsilon1) [below=of delta] {};
\node[ghost] (epsilon2) [below=of gamma] {};
\node[w2] at ($(epsilon1)!1/2!(epsilon2)$) {$\twid{\epsilon}_\Ho$};
\inedges{delta}{epsilon1}{\twid{r}_\Ho}
\inedge{gamma}{epsilon2}{\twid{\ell}_\Ho}
\inedge{delta}{gamma}{Fs}
\outabove{delta}{r_\Ho}
\outabove{gamma}{\ell_\Ho}
\outright{gamma}{t}
\outleft{delta}{t}
\end{doublestring}
=
\begin{doublestring}
\node[ghost] (ep1) {};
\node[ghost] (ep2) [right=of ep1] {};
\node[ghost] (epout1) [below=of ep1] {};
\node[ghost] (epout2) [below=of ep2] {};
\node[w2] at ($(ep1)!1/2!(ep2)$) {$\epsilon_\Ho$};
\throughhoriz{epout1}{epout2}{t}
\outabove{ep1}{r_\Ho}
\outabove{ep2}{\ell_\Ho}
\end{doublestring}
\end{equation}
Third and finally, $\gamma_\Ve$ and $\gamma_\Ho$ are compatible:
\begin{equation}
\label{eq:dbl-coh-1vc}
\begin{doublestring}
\node[ghost] (gammai1) {};
\node[ghost] (gammai2) [below=of gammai1] {};
\node[h2] at ($(gammai1)!1/2!(gammai2)$) {\makebox[0mm]{$\gamma_\Ve^{-1}$}};
\node[basic] (sigmat) [left=of gammai2] {$\twid{\sigma}$};
\node[ghost] (sigmatout) [above=of sigmat] {};
\inedges{sigmat}{gammai2}{\twid{\ell}_\Ve}
\outright{gammai1}{\ell_\Ve}
\outright{gammai2}{t}
\outbelow{sigmat}{\twid{\ell}_\Ho}
\outleftx{gammai1}{sigmatout}{Fs}
\end{doublestring}
=
\begin{doublestring}
\node[basic] (gamma) {$\gamma_\Ho$};
\node[basic] (sigma) [above=of gamma] {$\sigma$};
\inedge{sigma}{gamma}{\ell_\Ho}
\outright{sigma}{\ell_\Ve}
\outright{gamma}{t}
\outbelow{gamma}{\twid{\ell}_\Ho}
\outleft{gamma}{Fs}
\end{doublestring}
\qquad
\begin{doublestring}
\node[ghost] (gammai1) {};
\node[ghost] (gammai2) [below=of gammai1] {};
\node[h2] at ($(gammai1)!1/2!(gammai2)$) {\makebox[0mm]{$\gamma_\Ve^{-1}$}};
\node[basic] (tau) [right=of gammai1] {$\tau$};
\node[ghost] (tauout) [below=of tau] {};
\inedge{gammai1}{tau}{\ell_\Ve}
\outabove{tau}{\ell_\Ho}
\outrightx{gammai2}{tauout}{t}
\outleft{gammai2}{\twid{\ell}_\Ve}
\outleft{gammai1}{Fs}
\end{doublestring}
=
\begin{doublestring}
\node[basic] (gamma) {$\gamma_\Ho$};
\node[basic] (taut) [below=of gamma] {$\twid{\tau}$};
\inedge{gamma}{taut}{\twid{\ell}_\Ho}
\outabove{gamma}{\ell_\Ho}
\outright{gamma}{t}
\outleft{taut}{\twid{\ell}_\Ve}
\outleft{gamma}{Fs}
\end{doublestring}
\end{equation}
where $\gamma_\Ve^{-1}$ denotes the horizontal inverse of $\gamma_\Ve$.

A $2$-cell 
\[
\begin{tikzcd}[column sep=large]
(a, b, \ell, r, \eta, \epsilon, \sigma, \tau) 
\ar[r, "{(f, g, \lambda, \rho)}"] \ar[d, "{(s, t, \gamma, \delta)}"']
\ar[rd, phantom, "{(\alpha, \beta)}" description]       &
(a', b', \ldots)
\ar[d, "{(s', t', \gamma', \delta')}"]  \\
(\twid{a}, \twid{b}, \ldots)
\ar[r, "{(\twid{f}, \twid{g}, \twid{\lambda}, \twid{\rho})}"']   &
(\twid{a}', \twid{b}', \ldots)
\end{tikzcd}
\]
in $\C$ consists of $2$-cells
\[
\begin{tikzcd}
a \ar[r, "f"] \ar[d, "s"']
\ar[rd, phantom, "\alpha" description]  &
a' \ar[d, "s'"]         \\
\twid{a} \ar[r, "\twid{f}"']    &
\twid{a}'
\end{tikzcd}
\qquad
\begin{tikzcd}
b \ar[r, "g"] \ar[d, "t"']
\ar[rd, phantom, "\beta" description]   &
b' \ar[d, "t'"] \\
\twid{b} \ar[r, "\twid{g}"']    &
\twid{b}'
\end{tikzcd}
\]
in $\A$ and $\B$ respectively, satisfying the coherence axiom
\begin{equation}
\label{eq:dbl-coh-2}
\begin{doublestring}
\node[basic] (Falpha) {$F\alpha$};
\node[basic] (gamma') [right=of Falpha] {$\gamma'_\Ho$};
\node[ghost] (lambdat1) [below=of Falpha] {};
\node[ghost] (lambdat2) [below=of gamma'] {};
\node[w2] at ($(lambdat1)!1/2!(lambdat2)$) {\raisebox{-1mm}{$\twid{\lambda}_\Ho$}};
\inedges{Falpha}{lambdat1}{F\twid{f}}
\inedge{gamma'}{lambdat2}{\twid{\ell}'_\Ho}
\inedge{Falpha}{gamma'}{Fs'}
\outabove{Falpha}{Ff} 
\outabove{gamma'}{\ell'_\Ho}
\outright{gamma'}{t'}
\outbelow{lambdat2}{\twid{g}}
\outbelow{lambdat1}{\twid{\ell}_\Ho}
\outleft{Falpha}{Fs}
\end{doublestring}
=
\begin{doublestring}
\node[basic] (gamma) {$\gamma_\Ho$};
\node[basic] (beta) [right=of gamma] {$\beta$};
\node[ghost] (lambda1) [above=of gamma] {};
\node[ghost] (lambda2) [above=of beta] {};
\node[w2] at ($(lambda1)!1/2!(lambda2)$) {$\lambda_\Ho$};
\inedges{lambda1}{gamma}{\ell_\Ho}
\inedge{lambda2}{beta}{g}
\inedges{gamma}{beta}{t}
\outabove{lambda1}{Ff}
\outabove{lambda2}{\ell'_\Ho}
\outright{beta}{t'}
\outbelow{beta}{\twid{g}}
\outbelow{gamma}{\twid{\ell}_\Ho}
\outleft{gamma}{Fs}
\end{doublestring}
\end{equation}
Equation~\eqref{eq:dbl-coh-2} is just one of four similar coherence axioms,
but it is shown in Appendix~\ref{app:coh} that they are all equivalent.

Composition and identities in $\C$ are defined in the obvious way, making
use of the coherence isomorphisms of the double pseudofunctor $F$. In this
way, $\C$ is a double category, and the first and second projections
\[
\A \otby{P} \C \toby{Q} \B
\]
are strict double functors.

It remains to show that $P$ and $Q$ are surjective equivalences. As we
would expect by now, only the proof for $Q$ requires the hypothesis that
$F$ is a gregarious double equivalence.

$P$ is surjective on objects: given $a \in \A$, taking the identity
horizontal and vertical adjoint equivalences on $Fa$, and the identity
companion relation between the left adjoints, gives an object of $\C$ whose
image under $P$ is equal to $a$.

$P$ is horizontally full: take objects $(a, b, \ell, r, \eta, \epsilon,
\sigma, \tau)$ and $(a', b', \ldots)$ of $\C$, and a horizontal $1$-cell $a
\toby{f} a'$ in $\A$. Define
\[
g = \Bigl( b \toby{r_\Ho} Fa \toby{Ff} Fa' \toby{\ell'_\Ho} b \Bigr),
\]
then define
\[
\lambda_\Ho
=
\begin{doublestring}
\node[ghost] (eta1) {};
\node[ghost] (eta2) [right=of eta1] {};
\node[ghost] (etaout1) [right=of eta2] {};
\node[ghost] (etaout2) [right=of etaout1] {};
\node[w2] at ($(eta1)!1/2!(eta2)$) {$\eta_\Ho$};
\node[ghost] (e1) [below=of eta2] {};
\node[ghost] (e2) [right=of e1] {};
\node[ghost] (e3) [right=of e2] {};
\node[ghost] (eout) [left=of e1] {};
\node[w3] at (e2) {$=$};
\inedge{eta2}{e1}{r_\Ho}
\outabovex{e2}{etaout1}{Ff}
\outabovex{e3}{etaout2}{\ell'_\Ho}
\outbelow{e2}{g}
\outbelowx{eta1}{eout}{\ell_\Ho}
\end{doublestring}
\qquad
\rho_\Ho
=
\begin{doublestring}
\node[ghost] (eta'1) {};
\node[ghost] (eta'2) [right=of eta'1] {};
\node[ghost] (eta'out2) [left=of eta'1] {};
\node[ghost] (eta'out1) [left=of eta'out2] {};
\node[w2] at ($(eta'1)!1/2!(eta'2)$) {$\eta'_\Ho$};
\node[ghost] (e3) [below=of eta'1] {};
\node[ghost] (e2) [left=of e3] {};
\node[ghost] (e1) [left=of e2] {};
\node[ghost] (eout) [right=of e3] {};
\node[w3] at (e2) {$=$};
\inedge{eta'1}{e3}{\ell'_\Ho}
\outabovex{e1}{eta'out1}{r_\Ho}
\outabovex{e2}{eta'out2}{Ff}
\outbelow{e2}{g}
\outbelowx{eta2}{eout}{r'_\Ho}
\end{doublestring},
\]
which are vertically invertible since $\eta_\Ho$ and $\eta'_\Ho$ are. Also put
\[
\lambda_\Ve
=
\begin{doublestring}
\node[ghost] (eta1) {};
\node[ghost] (eta2) [right=of eta1] {};
\node[ghost] (etaout) [right=of eta2] {};
\node[w2] at ($(eta1)!1/2!(eta2)$) {$\eta_\Ho$};
\node[basic] (tau) [below=of eta1] {$\tau$};
\node[ghost] (etae) [below=of eta2] {};
\node[ghost] (e1) [below=of etae] {};
\node[ghost] (e2) [right=of e1] {};
\node[ghost] (e3) [right=of e2] {};
\node[w3] at (e2) {$=$};
\node[basic] (sigma') [above=of e3] {$\sigma'$};
\inedge{eta1}{tau}{\ell_\Ho}
\inedge{eta2}{e1}{r_\Ho}
\inedges{sigma'}{e3}{\ell'_\Ho}
\outabovex{e2}{etaout}{Ff}
\outright{sigma'}{\ell'_\Ve}
\outbelow{e2}{g}
\outleft{tau}{\ell_\Ve}
\end{doublestring}
\qquad
\rho_\Ve
=
\begin{doublestring}
\node[ghost] (e1) {};
\node[ghost] (e2) [right=of e1] {};
\node[ghost] (e3) [right=of e2] {};
\node[ghost] (e23) at ($(e2)!1/2!(e3)$) {};
\node[ghost] (e4) [right=of e3] {};
\node[w4] at (e23) {$=$};
\node[basic] (sigmabar) [below=of e1] {$\ovln{\sigma}$};
\node[basic] (tau') [below=of e4] {$\tau'$};
\node[ghost] (tau'out) [below=of tau'] {};
\node[ghost] (eta'i1) [left=of tau'] {};
\node[ghost] (eta'i2) [below=of eta'i1] {};
\node[ghost] (eta'iout1) [left=of eta'i2] {};
\node[ghost] (eta'iout2) [right=of eta'i2] {};
\node[h2] at ($(eta'i1)!1/2!(eta'i2)$) {\makebox[0mm]{$\!\!{\eta'}^{-1}_\Ve$}};
\inedge{e1}{sigmabar}{r_\Ho}
\inedge{e4}{tau'}{\ell'_\Ho}
\inedge{eta'i1}{tau'}{\ell'_\Ve}
\outabove{e23}{g}
\outrightx{eta'i2}{eta'iout2}{r'_\Ve}
\outbelowx{e2}{eta'iout1}{Ff}
\outleft{sigmabar}{r_\Ve}
\end{doublestring},
\]
where, in the last equation, $\ovln{\sigma}$ is as defined in
equations~\eqref{eq:dbl-mates} and ${\eta'}_\Ve^{-1}$ denotes the horizontal
inverse of $\eta'_\Ve$. It is routine to check
the coherence axioms, so that $(f, g, \lambda, \rho)$ is a horizontal
$1$-cell in $\C$ whose image under $P$ is equal to $f$. 

$P$ is vertically full, similarly.

$P$ is full and faithful on $2$-cells: take $0$- and $1$-cells 
\begin{equation}
\label{eq:dbl-2-shape}
\begin{tikzcd}[column sep=large]
(a, b, \ell, r, \eta, \epsilon, \sigma, \tau) 
\ar[r, "{(f, g, \lambda, \rho)}"] \ar[d, "{(s, t, \gamma, \delta)}"']   &
(a', b', \ldots)
\ar[d, "{(s', t', \gamma', \delta')}"]  \\
(\twid{a}, \twid{b}, \ldots)
\ar[r, "{(\twid{f}, \twid{g}, \twid{\lambda}, \twid{\rho})}"']   &
(\twid{a}', \twid{b}', \ldots)
\end{tikzcd}
\end{equation}
in $\C$ and a $2$-cell
\[
\begin{tikzcd}
a \ar[r, "f"] \ar[d, "s"']
\ar[rd, phantom, "\alpha" description]  &
a' \ar[d, "s'"]         \\
\twid{a} \ar[r, "\twid{f}"']    &
\twid{a}'
\end{tikzcd}
\]
in $\A$. We must prove that there is a unique $2$-cell
\[
\begin{tikzcd}
b \ar[r, "g"] \ar[d, "t"']
\ar[rd, phantom, "\beta" description]   &
b' \ar[d, "t'"] \\
\twid{b} \ar[r, "\twid{g}"']    &
\twid{b}'
\end{tikzcd}
\]
in $\B$ satisfying equation~\eqref{eq:dbl-coh-2}.

First we prove a general fact: for any $2$-cell of the form
\[
\begin{doublestring}
\node[ghost] (theta1) {};
\node[ghost] (theta2) [right=of theta1] {};
\node[w2] at ($(theta1)!1/2!(theta2)$) {$\theta$};
\outabove{theta1}{\ell_\Ho}
\outabove{theta2}{}
\outright{theta2}{}
\outbelow{theta2}{}
\outbelow{theta1}{\twid{\ell}_\Ho}
\outleft{theta1}{Fs}
\end{doublestring}
\]
the equation
\begin{equation}
\label{eq:dbl-canc}
\begin{doublestring}
\node[basic] (gamma) {$\gamma_\Ho$};
\node[basic] (xi) [right=of gamma] {$\xi$};
\inedge{gamma}{xi}{}
\outabove{gamma}{}
\outabove{xi}{}
\outright{xi}{}
\outbelow{xi}{}
\outbelow{gamma}{}
\outleft{gamma}{}
\end{doublestring}
=
\begin{doublestring}
\node[ghost] (theta1) {};
\node[ghost] (theta2) [right=of theta1] {};
\node[w2] at ($(theta1)!1/2!(theta2)$) {$\theta$};
\outabove{theta1}{}
\outabove{theta2}{}
\outright{theta2}{}
\outbelow{theta2}{}
\outbelow{theta1}{}
\outleft{theta1}{}
\end{doublestring}
\end{equation}
has exactly one solution $\xi$. Indeed,
\begin{align}
\label{eq:dc-1}
\eqref{eq:dbl-canc}
&
\implies
\begin{doublestring}
\node[basic] (delta) {$\delta_\Ho$};
\node[basic] (gamma) [right=of delta] {$\gamma_\Ho$};
\node[basic] (xi) [right=of gamma] {$\xi$};
\inedge{delta}{gamma}{}
\inedge{gamma}{xi}{}
\outabove{delta}{}
\outabove{gamma}{}
\outabove{xi}{}
\outright{xi}{}
\outbelow{xi}{}
\outbelow{gamma}{}
\outbelow{delta}{}
\outleft{delta}{}
\end{doublestring}
=
\begin{doublestring}
\node[basic] (delta) {$\delta_\Ho$};
\node[ghost] (theta1) [right=of delta] {};
\node[ghost] (theta2) [right=of theta1] {};
\node[w2] at ($(theta1)!1/2!(theta2)$) {$\theta$};
\inedge{delta}{theta1}{}
\outabove{delta}{}
\outabove{theta1}{}
\outabove{theta2}{}
\outright{theta2}{}
\outbelow{theta2}{}
\outbelow{theta1}{}
\outbelow{delta}{}
\outleft{delta}{}
\end{doublestring}
\\
\label{eq:dc-2}
&
\iff
\begin{doublestring}
\node[basic] (delta) {$\delta_\Ho$};
\node[basic] (gamma) [right=of delta] {$\gamma_\Ho$};
\node[basic] (xi) [right=of gamma] {$\xi$};
\node[ghost] (ept1) [below=of delta] {};
\node[ghost] (ept2) [right=of ept1] {};
\node[ghost] (eptout) [right=of ept2] {};
\node[w2] at ($(ept1)!1/2!(ept2)$) {$\twid{\epsilon}_\Ho$};
\inedge{delta}{gamma}{}
\inedge{gamma}{xi}{}
\inedge{delta}{ept1}{}
\inedge{gamma}{ept2}{}
\outabove{delta}{}
\outabove{gamma}{}
\outabove{xi}{}
\outright{xi}{}
\outbelow{xi}{}
\outleft{delta}{}
\end{doublestring}
=
\begin{doublestring}
\node[basic] (delta) {$\delta_\Ho$};
\node[ghost] (theta1) [right=of delta] {};
\node[ghost] (theta2) [right=of theta1] {};
\node[w2] at ($(theta1)!1/2!(theta2)$) {$\theta$};
\node[ghost] (ept1) [below=of delta] {};
\node[ghost] (ept2) [right=of ept1] {};
\node[ghost] (eptout) [right=of ept2] {};
\node[w2] at ($(ept1)!1/2!(ept2)$) {$\twid{\epsilon}_\Ho$};
\inedge{delta}{theta1}{}
\inedge{delta}{ept1}{}
\inedge{theta1}{ept2}{}
\outabove{delta}{}
\outabove{theta1}{}
\outabove{theta2}{}
\outright{theta2}{}
\outbelow{theta2}{}
\outleft{delta}{}
\end{doublestring}
\\
\label{eq:dc-3}
&
\iff
\begin{doublestring}
\node[ghost] (ep1) {};
\node[ghost] (ep2) [right=of ep1] {};
\node[ghost] (epout1) [below=of ep1] {};
\node[ghost] (epout2) [right=of ep2] {};
\node[w2] at ($(ep1)!1/2!(ep2)$) {$\epsilon_\Ho$};
\node[basic] (xi) [below=of epout2] {$\xi$};
\outabove{ep1}{}
\outabove{ep2}{}
\outabovex{xi}{epout2}{}
\outright{xi}{}
\outbelow{xi}{}
\outleftx{xi}{epout1}{}
\end{doublestring}
=
\begin{doublestring}
\node[basic] (delta) {$\delta_\Ho$};
\node[ghost] (theta1) [right=of delta] {};
\node[ghost] (theta2) [right=of theta1] {};
\node[w2] at ($(theta1)!1/2!(theta2)$) {$\theta$};
\node[ghost] (ept1) [below=of delta] {};
\node[ghost] (ept2) [right=of ept1] {};
\node[ghost] (eptout) [right=of ept2] {};
\node[w2] at ($(ept1)!1/2!(ept2)$) {$\twid{\epsilon}_\Ho$};
\inedge{delta}{theta1}{}
\inedge{delta}{ept1}{}
\inedge{theta1}{ept2}{}
\outabove{delta}{}
\outabove{theta1}{}
\outabove{theta2}{}
\outright{theta2}{}
\outbelow{theta2}{}
\outleft{delta}{}
\end{doublestring}
\\
\label{eq:dc-4}
&
\iff
\begin{doublestring}
\node[basic] (xi) {$\xi$};
\outabove{xi}{}
\outright{xi}{}
\outbelow{xi}{}
\outleft{xi}{}
\end{doublestring}
=
\begin{doublestring}
\node[ghost] (epi1) {};
\node[ghost] (epi2) [right=of epi1] {};
\node[ghost] (epiout) [right=of epi2] {};
\node[w2] at ($(epi1)!1/2!(epi2)$) {$\epsilon_\Ho^{-1}$};
\node[basic] (delta) [below=of epi1] {$\delta_\Ho$};
\node[ghost] (theta1) [right=of delta] {};
\node[ghost] (theta2) [right=of theta1] {};
\node[w2] at ($(theta1)!1/2!(theta2)$) {$\theta$};
\node[ghost] (ept1) [below=of delta] {};
\node[ghost] (ept2) [right=of ept1] {};
\node[ghost] (eptout) [right=of ept2] {};
\node[w2] at ($(ept1)!1/2!(ept2)$) {$\twid{\epsilon}_\Ho$};
\inedge{epi1}{delta}{}
\inedge{epi2}{theta1}{}
\inedge{delta}{ept1}{}
\inedge{theta1}{ept2}{}
\outabovex{theta2}{epiout}{}
\outright{theta2}{}
\outbelowx{theta2}{eptout}{}
\outleft{delta}{}
\end{doublestring},
\end{align}
where the double implication~\eqref{eq:dc-2} holds because
$\twid{\epsilon}_\Ho$ is vertically invertible, \eqref{eq:dc-3} follows
from the second equation of~\eqref{eq:dbl-coh-1vb}, and~\eqref{eq:dc-4}
holds because $\epsilon_\Ho$ is vertically invertible. This proves
that~\eqref{eq:dbl-canc} has at most one solution $\xi$, or in other words,
that $\gamma_\Ho$ is left cancellable. A similar argument show that
$\delta_\Ho$ is left cancellable. The implication in~\eqref{eq:dc-1}
therefore goes in both directions. Hence, reading upwards, the $2$-cell
$\xi$ defined by~\eqref{eq:dc-4} is indeed a solution
of~\eqref{eq:dbl-canc}.

It now follows that $P$ is full and faithful on $2$-cells:
for~\eqref{eq:dbl-coh-2} is equivalent to 
\[
\begin{doublestring}
\node[basic] (gamma) {$\gamma_\Ho$};
\node[basic] (beta) [right=of gamma] {$\beta$};
\inedges{gamma}{beta}{}
\outabove{gamma}{}
\outabove{beta}{}
\outright{beta}{}
\outbelow{beta}{}
\outbelow{gamma}{}
\outleft{gamma}{}
\end{doublestring}
=
\begin{doublestring}
\node[basic] (Falpha) {$F\alpha$};
\node[basic] (gamma') [right=of Falpha] {$\gamma'_\Ho$};
\node[ghost] (lambdat1) [below=of Falpha] {};
\node[ghost] (lambdat2) [below=of gamma'] {};
\node[w2] at ($(lambdat1)!1/2!(lambdat2)$) {\raisebox{-1mm}{$\twid{\lambda}_\Ho$}};
\node[ghost] (lambdai1) [above=of Falpha] {};
\node[ghost] (lambdai2) [above=of gamma'] {};
\node[w2] at ($(lambdai1)!1/2!(lambdai2)$) {$\lambda_\Ho^{-1}$};
\inedge{Falpha}{lambdat1}{}
\inedge{gamma'}{lambdat2}{}
\inedge{Falpha}{gamma'}{}
\inedge{lambdai1}{Falpha}{}
\inedge{lambdai2}{gamma'}{}
\outabove{lambdai1}{}
\outabove{lambdai2}{}
\outright{gamma'}{}
\outbelow{lambdat2}{}
\outbelow{lambdat1}{}
\outleft{Falpha}{}
\end{doublestring}
\]
which by the general fact just proved has a unique solution $\beta$.

This completes the proof that $P$ is a surjective equivalence. We now prove
that $Q$ is too.

$Q$ is surjective on objects: let $b \in \B$. Since $F$ is surjective on
objects up to gregarious equivalence, we can choose an object $a \in \A$,
horizontal and vertical equivalences 
\[
\begin{tikzcd}
Fa \ar[r, "\ell_\Ho"] & b
\end{tikzcd}
\qquad
\begin{tikzcd}
Fa \ar[d, "\ell_\Ve"] \\ b
\end{tikzcd}
\]
in $\B$, and binding $2$-cells $\sigma$ and $\tau$ making $\ell_\Ho$ and
$\ell_\Ve$ into companions. By a standard $2$-categorical result,
$\ell_\Ho$ can be extended to an adjoint equivalence $(\ell_\Ho, r_\Ho,
\eta_\Ho, \epsilon_\Ho)$ in the horizontal $2$-category of $\B$, and
$\ell_\Ve$ similarly. Then $(a, b, \ell, r, \eta, \epsilon, \sigma, \tau)$
is an object of $\C$, and its image under $Q$ is equal to $b$.

$Q$ is horizontally full: take $0$-cells $(a, b, \ell, r, \eta, \epsilon,
\sigma, \tau)$ and $(a', b', \ldots)$ of $\C$, and a horizontal $1$-cell $b
\toby{g} b'$ in $\B$. Since $F$ is horizontally essentially full, there
exist a horizontal $1$-cell $a \toby{f} a'$ in $\A$ and a vertically
invertible $2$-cell
\[
\begin{tikzcd}
Fa \ar[rrr, "Ff"] \ar[d, equal]
\ar[rrrd, phantom, "\chi" description]   &
&
&
Fa' \ar[d, equal]       \\
Fa \ar[r, "\ell_\Ho"']  &
b \ar[r, "g"']  &
b' \ar[r, "r'_\Ho"']    &
Fa'
\end{tikzcd}
\]
in $\B$. Define $2$-cells
\[
\lambda_\Ho
=
\begin{doublestring}
\node[ghost] (chi1) {};
\node[ghost] (chi2) [right=of chi1] {};
\node[ghost] (chi3) [right=of chi2] {};
\node[ghost] (chiout) [right=of chi3] {};
\node[w3] at (chi2) {$\chi$};
\node[ghost] (ep'1) [below=of chi3] {};
\node[ghost] (ep'2) [right=of ep'1] {};
\node[ghost] (ep'out2) [left=of ep'1] {};
\node[ghost] (ep'out1) [left=of ep'out2] {};
\node[w2] at ($(ep'1)!1/2!(ep'2)$) {$\epsilon'_\Ho$};
\inedge{chi3}{ep'1}{r'_\Ho}
\outabove{chi2}{Ff}
\outabovex{ep'2}{chiout}{\ell'_\Ho}
\outbelowx{chi2}{ep'out2}{g}
\outbelowx{chi1}{ep'out1}{\ell_\Ho}
\end{doublestring}
\qquad
\rho_\Ho
=
\begin{doublestring}
\node[ghost] (chi1) {};
\node[ghost] (chi2) [right=of chi1] {};
\node[ghost] (chi3) [right=of chi2] {};
\node[ghost] (chiout) [left=of chi1] {};
\node[w3] at (chi2) {$\chi$};
\node[ghost] (ep2) [below=of chi1] {};
\node[ghost] (ep1) [left=of ep2] {};
\node[ghost] (epout1) [right=of ep2] {};
\node[ghost] (epout2) [right=of epout1] {};
\node[w2] at ($(ep1)!1/2!(ep2)$) {$\epsilon_\Ho$};
\inedge{chi1}{ep2}{\ell_\Ho}
\outabovex{ep1}{chiout}{r_\Ho}
\outabove{chi2}{Ff}
\outbelowx{chi3}{epout2}{r'_\Ho}
\outbelowx{chi2}{epout1}{g}
\end{doublestring}
\]
in $\B$, which are vertically invertible since $\chi$, $\epsilon_\Ho$ and
$\epsilon'_\Ho$ are. Also put
\[
\lambda_\Ve
=
\begin{doublestring}
\node[ghost] (chi1) {};
\node[ghost] (chi2) [right=of chi1] {};
\node[ghost] (chi3) [right=of chi2] {};
\node[w3] at (chi2) {$\chi$};
\node[basic] (tau) [below=of chi1] {$\tau$};
\node[ghost] (tauout) [right=of tau] {};
\node[ghost] (ep'1) [below=of chi3] {};
\node[ghost] (ep'2) [right=of ep'1] {};
\node[w2] at ($(ep'1)!1/2!(ep'2)$) {$\epsilon'_\Ho$};
\node[basic] (sigma') [above=of ep'2] {$\sigma'$};
\inedge{chi1}{tau}{\ell_\Ho}
\inedges{chi3}{ep'1}{r'_\Ho}
\inedges{sigma'}{ep'2}{\ell'_\Ho}
\outabove{chi2}{Ff}
\outright{sigma'}{\ell'_\Ve}
\outbelowx{chi2}{tauout}{g}
\outleft{tau}{\ell_\Ve}
\end{doublestring}
\qquad
\rho_\Ve
=
\begin{doublestring}
\node[ghost] (chii1) {};
\node[ghost] (chii2) [right=of chii1] {};
\node[ghost] (chii3) [right=of chii2] {};
\node[ghost] (chii23) at ($(chii2)!1/2!(chii3)$) {};
\node[ghost] (chii4) [right=of chii3] {};
\node[w4] at (chii23) {$\chi^{-1}$};
\node[basic] (sigma) [above=of chii1] {$\sigma$};
\node[ghost] (sigmaout) [above=of sigma] {};
\node[basic] (tau'bar) [above=of chii4] {$\ovln{\tau'}$};
\node[ghost] (epi2) [right=of sigma] {};
\node[ghost] (epi1) [above=of epi2] {};
\node[ghost] (epiout2) [right=of epi2] {};
\node[ghost] (epiout1) [right=of epi1] {};
\node[h2] at ($(epi1)!1/2!(epi2)$) {\makebox[0mm]{$\epsilon^{-1}_\Ve$}};
\inedge{sigma}{chii1}{\ell_\Ho}
\inedge{sigma}{epi2}{\ell_\Ve}
\inedges{tau'bar}{chii4}{r'_\Ho}
\outabovex{chii3}{epiout1}{g}
\outright{tau'bar}{r'_\Ve}
\outbelow{chii23}{Ff}
\outleftx{epi1}{sigmaout}{r_\Ve}
\end{doublestring}
\]
where $\ovln{\tau'}$ is defined as in equations~\eqref{eq:dbl-mates}. A
series of routine checks proves that $(f, g, \lambda, \rho)$ is indeed a
horizontal $1$-cell $(a, b, \ldots) \to (a', b', \ldots)$ in $\C$, and
plainly its image under $Q$ is equal to $g$.

$Q$ is vertically full, similarly.

$Q$ is full and faithful on $2$-cells: take $0$- and
$1$-cells~\eqref{eq:dbl-2-shape} in $\C$ and a $2$-cell
\[
\begin{tikzcd}
b \ar[r, "g"] \ar[d, "t"']
\ar[rd, phantom, "\beta" description]   &
b' \ar[d, "t'"] \\
\twid{b} \ar[r, "\twid{g}"']    &
\twid{b}'
\end{tikzcd}
\]
in $\B$. We must prove that there is a unique $2$-cell
\[
\begin{tikzcd}
a \ar[r, "f"] \ar[d, "s"']
\ar[rd, phantom, "\alpha" description]  &
a' \ar[d, "s'"]         \\
\twid{a} \ar[r, "\twid{f}"']    &
\twid{a}'
\end{tikzcd}
\]
in $\A$ satisfying equation~\eqref{eq:dbl-coh-2}. This is proved by the
same argument as was used to show that $P$ is full and faithful on
$2$-cells, together with the hypothesis that $F$ is full and faithful on
$2$-cells. 

This completes the proof that $Q$ is a surjective equivalence, and,
therefore, of Theorem~\ref{thm:D}.
\end{pfof}

\begin{cor}
\label{cor:sym}
Gregarious double equivalence is a symmetric relation on double categories.
\end{cor}

That is, for double categories $\A$ and $\B$, if there exists a gregarious
double equivalence $\A \to \B$ then there exists a gregarious double
equivalence $\B \to \A$. This can also be proved directly, by a more
complicated version of the proof of Lemma~\ref{lemma:inverse}.

\paragraph*{Acknowledgements} I thank Nathanael Arkor, John Bourke,
Alexander Campbell, Bryce Clarke, St\'ephane Desarzens, Simon Henry, Emily
Riehl, Maru Sarazola, Mike Shulman and Carlos Simpson for helpful comments.

\appendix

\section{Appendix: Coherence conditions}
\label{app:coh}

Here we prove the claim made after equation~\eqref{eq:dbl-coh-2} in the
proof of Theorem~\ref{thm:D}. There, we were defining the $2$-cells in the
double category $\C$. The claim was that equation~\eqref{eq:dbl-coh-2} was
just one of four similar coherence equations, all of which are equivalent.

First we state those four equations, including~\eqref{eq:dbl-coh-2} itself
for convenience:
\begin{align}
\label{eq:dbl-coh-2a}
\begin{doublestring}
\node[basic] (Falpha) {$F\alpha$};
\node[basic] (gamma') [right=of Falpha] {$\gamma'_\Ho$};
\node[ghost] (lambdat1) [below=of Falpha] {};
\node[ghost] (lambdat2) [below=of gamma'] {};
\node[w2] at ($(lambdat1)!1/2!(lambdat2)$) {\raisebox{-1mm}{$\twid{\lambda}_\Ho$}};
\inedges{Falpha}{lambdat1}{}%{F\twid{f}}
\inedge{gamma'}{lambdat2}{}%{\twid{\ell}'_\Ho}
\inedge{Falpha}{gamma'}{}%{Fs'}
\outabove{Falpha}{}%{Ff} 
\outabove{gamma'}{}%{\ell'_\Ho}
\outright{gamma'}{}%{t'}
\outbelow{lambdat2}{}%{\twid{g}}
\outbelow{lambdat1}{}%{\twid{\ell}_\Ho}
\outleft{Falpha}{}%{Fs}
\end{doublestring}
&
=
\begin{doublestring}
\node[basic] (gamma) {$\gamma_\Ho$};
\node[basic] (beta) [right=of gamma] {$\beta$};
\node[ghost] (lambda1) [above=of gamma] {};
\node[ghost] (lambda2) [above=of beta] {};
\node[w2] at ($(lambda1)!1/2!(lambda2)$) {$\lambda_\Ho$};
\inedges{lambda1}{gamma}{}%{\ell_\Ho}
\inedge{lambda2}{beta}{}%{g}
\inedges{gamma}{beta}{}%{t}
\outabove{lambda1}{}%{Ff}
\outabove{lambda2}{}%{\ell'_\Ho}
\outright{beta}{}%{t'}
\outbelow{beta}{}%{\twid{g}}
\outbelow{gamma}{}%{\twid{\ell}_\Ho}
\outleft{gamma}{}%{Fs}
\end{doublestring}
\\
\label{eq:dbl-coh-2b}
\begin{doublestring}
\node[basic] (delta) {$\delta_\Ho$};
\node[basic] (Falpha) [right=of delta] {$F\alpha$};
\node[ghost] (rhot1) [below=of delta] {};
\node[ghost] (rhot2) [below=of Falpha] {};
\node[w2] at ($(rhot1)!1/2!(rhot2)$) {$\twid{\rho}_\Ho$};
\inedges{delta}{rhot1}{}%{r_H}
\inedge{Falpha}{rhot2}{}%{Ff}
\inedge{delta}{Falpha}{}%{Fs}
\outabove{delta}{}%{r_H}
\outabove{Falpha}{}%{Ff}
\outright{Falpha}{}%{Fs'}
\outbelow{rhot2}{}%{\twid{r}'_\Ho}
\outbelow{rhot1}{}%{\twid{g}}
\outleft{delta}{}%{t}
\end{doublestring}
&
=
\begin{doublestring}
\node[basic] (beta) {$\beta$};
\node[basic] (delta') [right=of beta] {$\delta'_\Ho$};
\node[ghost] (rho1) [above=of beta] {};
\node[ghost] (rho2) [above=of delta'] {};
\node[w2] at ($(rho1)!1/2!(rho2)$) {$\rho_\Ho$};
\inedges{rho1}{beta}{}%{g}
\inedge{rho2}{delta'}{}%{r'_\Ho}
\inedges{beta}{delta'}{}%{t'}
\outabove{rho1}{}%{r_\Ho}
\outabove{rho2}{}%{Ff}
\outright{delta'}{}%{Fs'}
\outbelow{delta'}{}%{\twid{r}'_\Ho}
\outbelow{beta}{}%{\twid{g}}
\outleft{beta}{}%{t}
\end{doublestring}
\\
\label{eq:dbl-coh-2c}
\begin{doublestring}
\node[basic] (lambda) {$\lambda_\Ve$};
\node[basic] (beta) [below=of lambda] {$\beta$};
\node[ghost] (gamma'1) [right=of lambda] {};
\node[ghost] (gamma'2) [right=of beta] {};
\node[h2] at ($(gamma'1)!1/2!(gamma'2)$) {$\gamma'_\Ve$};
\inedges{lambda}{beta}{}%{g}
\inedge{lambda}{gamma'1}{}%{\ell'_\Ve}
\inedges{beta}{gamma'2}{}%{t'}
\outabove{lambda}{}%{Ff}
\outright{gamma'1}{}%{Fs'}
\outright{gamma'2}{}%{\twid{\ell}'_\Ve}
\outbelow{beta}{}%{\twid{g}}
\outleft{beta}{}%{t}
\outleft{lambda}{}%{\ell_\Ve}
\end{doublestring}
&
=
\begin{doublestring}
\node[basic] (Falpha) {$F\alpha$};
\node[basic] (lambdat) [below=of Falpha] {\raisebox{-1mm}{$\twid{\lambda}_\Ve$}};
\node[ghost] (gamma1) [left=of Falpha] {};
\node[ghost] (gamma2) [left=of lambdat] {};
\node[h2] at ($(gamma1)!1/2!(gamma2)$) {$\gamma_\Ve$};
\inedge{Falpha}{lambdat}{}%{F\twid{f}}
\inedge{gamma1}{Falpha}{}%{Fs}
\inedges{gamma2}{lambdat}{}%{\twid{\ell}_\Ve}
\outabove{Falpha}{}%{Ff}
\outright{Falpha}{}%{Fs'}
\outright{lambdat}{}%{\twid{\ell}'_\Ve}
\outbelow{lambdat}{}%{\twid{g}}
\outleft{gamma2}{}%{t}
\outleft{gamma1}{}%{\ell_\Ve}
\end{doublestring}
\\
\label{eq:dbl-coh-2d}
\begin{doublestring}
\node[basic] (beta) {$\beta$};
\node[basic] (rhot) [below=of beta] {$\twid{\rho}_\Ve$};
\node[ghost] (delta'1) [right=of beta] {};
\node[ghost] (delta'2) [right=of rhot] {};
\node[h2] at ($(delta'1)!1/2!(delta'2)$) {$\delta'_\Ve$};
\inedges{beta}{rhot}{}%{\twid{g}}
\inedge{beta}{delta'1}{}%{t'}
\inedges{rhot}{delta'2}{}%{\twid{r}'_\Ve}
\outabove{beta}{}%{g}
\outright{delta'1}{}%{r'_\Ve}
\outright{delta'2}{}%{Fs'}
\outbelow{rhot}{}%{F\twid{f}}
\outleft{rhot}{}%{\twid{r}_\Ve}
\outleft{beta}{}%{t}
\end{doublestring}
&
=
\begin{doublestring}
\node[basic] (rho) {$\rho_\Ve$};
\node[basic] (Falpha) [below=of rho] {$F\alpha$};
\node[ghost] (delta1) [left=of rho] {};
\node[ghost] (delta2) [left=of Falpha] {};
\node[h2] at ($(delta1)!1/2!(delta2)$) {$\delta_\Ve$};
\inedge{rho}{Falpha}{}%{Ff}
\inedge{delta1}{rho}{}%{r_\Ve}
\inedges{delta2}{Falpha}{}%{Fs}
\outabove{rho}{}%{g}
\outright{rho}{}%{r'_\Ve}
\outright{Falpha}{}%{Fs'}
\outbelow{Falpha}{}%{F\twid{f}}
\outleft{delta2}{}%{\twid{r}_\Ve}
\outleft{delta1}{}%{t}
\end{doublestring}
\end{align}
By symmetry, in order to show that
\eqref{eq:dbl-coh-2a}--\eqref{eq:dbl-coh-2d} are equivalent, it suffices to
show that~\eqref{eq:dbl-coh-2a} implies all the others. Moreover, if we can
show that \eqref{eq:dbl-coh-2a}$\textimplies$\eqref{eq:dbl-coh-2b}, then it
follows by symmetry again that
\eqref{eq:dbl-coh-2c}$\textimplies$\eqref{eq:dbl-coh-2d}. Hence it suffices
to show that~\eqref{eq:dbl-coh-2a} implies both~\eqref{eq:dbl-coh-2b}
and~\eqref{eq:dbl-coh-2c}.

Assuming~\eqref{eq:dbl-coh-2a}, we first prove~\eqref{eq:dbl-coh-2b}:
\begin{align*}
\begin{doublestring}
\node[basic] (delta) {$\delta_\Ho$};
\node[basic] (Falpha) [right=of delta] {$F\alpha$};
\node[ghost] (rhot1) [below=of delta] {};
\node[ghost] (rhot2) [below=of Falpha] {};
\node[w2] at ($(rhot1)!1/2!(rhot2)$) {$\twid{\rho}_\Ho$};
\inedges{delta}{rhot1}{}%{r_H}
\inedge{Falpha}{rhot2}{}%{Ff}
\inedge{delta}{Falpha}{}%{Fs}
\outabove{delta}{}%{r_H}
\outabove{Falpha}{}%{Ff}
\outright{Falpha}{}%{Fs'}
\outbelow{rhot2}{}%{\twid{r}'_\Ho}
\outbelow{rhot1}{}%{\twid{g}}
\outleft{delta}{}%{t}
\end{doublestring}
&
=
\begin{doublestring}
\node[basic] (delta) {\makebox[0mm]{$\delta_\Ho$}};
\node[basic] (Falpha) [right=of delta] {\makebox[0mm]{$F\alpha$}};
\node[ghost] (rhot1) [below=of delta] {};
\node[ghost] (rhot2) [below=of Falpha] {};
\node[w2] at ($(rhot1)!1/2!(rhot2)$) {$\twid{\rho}_\Ho$};
\node[ghost] (ept'1) [below=of rhot2] {};
\node[ghost] (ept'2) [right=of ept'1] {};
\node[ghost] (ept'out1) [left=of ept'1] {};
\node[ghost] (ept'out2) [right=of ept'2] {};
\node[w2] at ($(ept'1)!1/2!(ept'2)$) {$\twid{\epsilon}'_\Ho$};
\node[ghost] (etat'1) [above=of ept'2] {};
\node[ghost] (etat'2) [right=of etat'1] {};
\node[ghost] (etat'out) [above=of etat'2] {};
\node[w2] at ($(etat'1)!1/2!(etat'2)$) {$\twid{\eta}'_\Ho$};
\inedge{delta}{rhot1}{}
\inedge{Falpha}{rhot2}{}
\inedge{delta}{Falpha}{}
\inedge{rhot2}{ept'1}{}
\inedge{etat'1}{ept'2}{}
\outabove{delta}{}
\outabove{Falpha}{}
\outrightx{Falpha}{etat'out}{}
\outbelowx{etat'2}{ept'out2}{}
\outbelowx{rhot1}{ept'out1}{}
\outleft{delta}{}
\end{doublestring}
=
\begin{doublestring}
\node[basic] (delta) {\makebox[0mm]{$\delta_\Ho$}};
\node[basic] (Falpha) [right=of delta] {\makebox[0mm]{$F\alpha$}};
\node[ghost] (lambdatpad) [below=of Falpha] {};
\node[ghost] (lambdat1) [below=of lambdatpad] {};
\node[ghost] (lambdat2) [right=of lambdat1] {};
\node[w2] at ($(lambdat1)!1/2!(lambdat2)$) {$\twid{\lambda}_\Ho$};
\node[ghost] (ept2) [below=of lambdat1] {};
\node[ghost] (ept1) [left=of ept2] {};
\node[w2] at ($(ept1)!1/2!(ept2)$) {$\twid{\epsilon}_\Ho$};
\node[ghost] (eptout1) [right=of ept2] {};
\node[ghost] (eptout2) [right=of eptout1] {};
\node[ghost] (etat'1) [above=of lambdat2] {};
\node[ghost] (etat'2) [right=of etat'1] {};
\node[ghost] (etat'out) [above=of etat'2] {};
\node[w2] at ($(etat'1)!1/2!(etat'2)$) {$\twid{\eta}'_\Ho$};
\inedge{delta}{Falpha}{}
\inedge{delta}{ept1}{}
\inedge{Falpha}{lambdat1}{}
\inedge{lambdat1}{ept2}{}
\inedge{etat'1}{lambdat2}{}
\outabove{delta}{}
\outabove{Falpha}{}
\outrightx{Falpha}{etat'out}{}
\outbelowx{etat'2}{eptout2}{}
\outbelowx{lambdat2}{eptout1}{}
\outleft{delta}{}
\end{doublestring}
\\
&
=
\begin{doublestring}
\node[basic] (delta) {\makebox[0mm]{$\delta_\Ho$}};
\node[basic] (Falpha) [right=of delta] {\makebox[0mm]{$F\alpha$}};
\node[ghost] (lambdat1) [below=of Falpha] {};
\node[ghost] (lambdat2) [right=of lambdat1] {};
\node[w2] at ($(lambdat1)!1/2!(lambdat2)$) {$\twid{\lambda}_\Ho$};
\node[ghost] (ept2) [below=of lambdat1] {};
\node[ghost] (ept1) [left=of ept2] {};
\node[w2] at ($(ept1)!1/2!(ept2)$) {$\twid{\epsilon}_\Ho$};
\node[ghost] (eptout1) [right=of ept2] {};
\node[ghost] (eptout2) [right=of eptout1] {};
\node[basic] (gamma') [above=of lambdat2] {\makebox[0mm]{$\gamma'_\Ho$}};
\node[basic] (delta') [right=of gamma'] {\makebox[0mm]{$\delta'_\Ho$}};
\node[ghost] (eta'1) [above=of gamma'] {};
\node[ghost] (eta'2) [right=of eta'1] {};
\node[ghost] (eta'out2) [left=of eta'1] {};
\node[ghost] (eta'out1) [left=of eta'out2] {};
\node[w2] at ($(eta'1)!1/2!(eta'2)$) {$\eta'_\Ho$};
\inedge{delta}{Falpha}{}
\inedge{delta}{ept1}{}
\inedge{Falpha}{lambdat1}{}
\inedge{lambdat1}{ept2}{}
\inedge{Falpha}{gamma'}{}
\inedge{gamma'}{lambdat2}{}
\inedge{gamma'}{delta'}{}
\inedge{eta'1}{gamma'}{}
\inedge{eta'2}{delta'}{}
\outabovex{delta}{eta'out1}{}
\outabovex{Falpha}{eta'out2}{}
\outright{delta'}{}
\outbelowx{delta'}{eptout2}{}
\outbelowx{lambdat2}{eptout1}{}
\outleft{delta}{}
\end{doublestring}
=
\begin{doublestring}
\node[basic] (delta) {\makebox[0mm]{$\delta_\Ho$}};
\node[basic] (gamma) [right=of delta] {\makebox[0mm]{$\gamma_\Ho$}};
\node[basic] (beta) [right=of gamma] {$\beta$};
\node[ghost] (lambda1) [above=of gamma] {};
\node[ghost] (lambda2) [right=of lambda1] {};
\node[w2] at ($(lambda1)!1/2!(lambda2)$) {$\lambda_\Ho$};
\node[ghost] (ept2) [below=of gamma] {};
\node[ghost] (ept1) [left=of ept2] {};
\node[w2] at ($(ept1)!1/2!(ept2)$) {$\twid{\epsilon}_\Ho$};
\node[ghost] (eptout1) [right=of ept2] {};
\node[ghost] (eptout2) [right=of eptout1] {};
\node[basic] (delta') [right=of beta] {\makebox[0mm]{$\delta'_\Ho$}};
\node[ghost] (eta'1) [above=of lambda2] {};
\node[ghost] (eta'2) [right=of eta'1] {};
\node[ghost] (eta'out2) [left=of eta'1] {};
\node[ghost] (eta'out1) [left=of eta'out2] {};
\node[w2] at ($(eta'1)!1/2!(eta'2)$) {$\eta'_\Ho$};
\inedge{delta}{gamma}{}
\inedge{delta}{ept1}{}
\inedge{gamma}{ept2}{}
\inedge{lambda1}{gamma}{}
\inedge{lambda2}{beta}{}
\inedge{gamma}{beta}{}
\inedge{beta}{delta'}{}
\inedge{eta'1}{lambda2}{}
\inedge{eta'2}{delta'}{}
\outabovex{delta}{eta'out1}{}
\outabovex{lambda1}{eta'out2}{}
\outright{delta'}{}
\outbelowx{delta'}{eptout2}{}
\outbelowx{beta}{eptout1}{}
\outleft{delta}{}
\end{doublestring}
\\
&
=
\begin{doublestring}
\node[ghost] (lambda1) {};
\node[ghost] (lambda2) [right=of lambda1] {};
\node[w2] at ($(lambda1)!1/2!(lambda2)$) {$\lambda_\Ho$};
\node[ghost] (ep2) [below=of lambda1] {};
\node[ghost] (ep1) [left=of ep2] {};
\node[ghost] (epout) [below=of ep1] {};
\node[w2] at ($(ep1)!1/2!(ep2)$) {$\epsilon_\Ho$};
\node[ghost] (eta'1) [above=of lambda2] {};
\node[ghost] (eta'2) [right=of eta'1] {};
\node[ghost] (eta'out2) [left=of eta'1] {};
\node[ghost] (eta'out1) [left=of eta'out2] {};
\node[w2] at ($(eta'1)!1/2!(eta'2)$) {$\eta'_\Ho$};
\node[ghost] (betapad) [below=of lambda2] {};
\node[basic] (beta) [below=of betapad] {$\beta$};
\node[basic] (delta') [right=of beta] {\makebox[0mm]{$\delta'_\Ho$}};
\inedge{lambda1}{ep2}{}
\inedge{eta'1}{lambda2}{}
\inedge{eta'2}{delta'}{}
\inedge{lambda2}{beta}{}
\inedge{beta}{delta'}{}
\outabovex{ep1}{eta'out1}{}
\outabovex{lambda1}{eta'out2}{}
\outright{delta'}{}
\outbelow{delta'}{}
\outbelow{beta}{}
\outleftx{beta}{epout}{}
\end{doublestring}
=
\begin{doublestring}
\node[ghost] (eta1) {};
\node[ghost] (eta2) [right=of eta1] {};
\node[ghost] (etaout1) [left=of eta1] {};
\node[ghost] (etaout2) [right=of eta2] {};
\node[w2] at ($(eta1)!1/2!(eta2)$) {$\eta_\Ho$};
\node[ghost] (rho1) [below=of eta2] {};
\node[ghost] (rho2) [right=of rho1] {};
\node[w2] at ($(rho1)!1/2!(rho2)$) {$\rho_\Ho$};
\node[ghost] (eppad) [below=of eta1] {};
\node[ghost] (ep2) [below=of eppad] {};
\node[ghost] (ep1) [left=of ep2] {};
\node[ghost] (epout) [below=of ep1] {};
\node[w2] at ($(ep1)!1/2!(ep2)$) {$\epsilon_\Ho$};
\node[ghost] (betapad) [below=of rho1] {};
\node[basic] (beta) [below=of betapad] {$\beta$};
\node[basic] (delta') [right=of beta] {\makebox[0mm]{$\delta'_\Ho$}};
\inedge{eta1}{ep2}{}
\inedge{eta2}{rho1}{}
\inedge{rho1}{beta}{}
\inedge{rho2}{delta'}{}
\inedge{beta}{delta'}{}
\outabovex{ep1}{etaout1}{}
\outabovex{rho2}{etaout2}{}
\outright{delta'}{}
\outbelow{delta'}{}
\outbelow{beta}{}
\outleftx{beta}{epout}{}
\end{doublestring}
\\
&
=
\begin{doublestring}
\node[basic] (beta) {$\beta$};
\node[basic] (delta') [right=of beta] {$\delta'_\Ho$};
\node[ghost] (rho1) [above=of beta] {};
\node[ghost] (rho2) [above=of delta'] {};
\node[w2] at ($(rho1)!1/2!(rho2)$) {$\rho_\Ho$};
\inedges{rho1}{beta}{}%{g}
\inedge{rho2}{delta'}{}%{r'_\Ho}
\inedges{beta}{delta'}{}%{t'}
\outabove{rho1}{}%{r_\Ho}
\outabove{rho2}{}%{Ff}
\outright{delta'}{}%{Fs'}
\outbelow{delta'}{}%{\twid{r}'_\Ho}
\outbelow{beta}{}%{\twid{g}}
\outleft{beta}{}%{t}
\end{doublestring},
\end{align*}
where the first and last equalities are by adjointness, the second and
second to last are by~\eqref{eq:dbl-coh-1ha}, the third and third to last
are by~\eqref{eq:dbl-coh-1vb}, and the middle equality is
by~\eqref{eq:dbl-coh-2a}.

Now assuming~\eqref{eq:dbl-coh-2a}, we prove~\eqref{eq:dbl-coh-2c}:
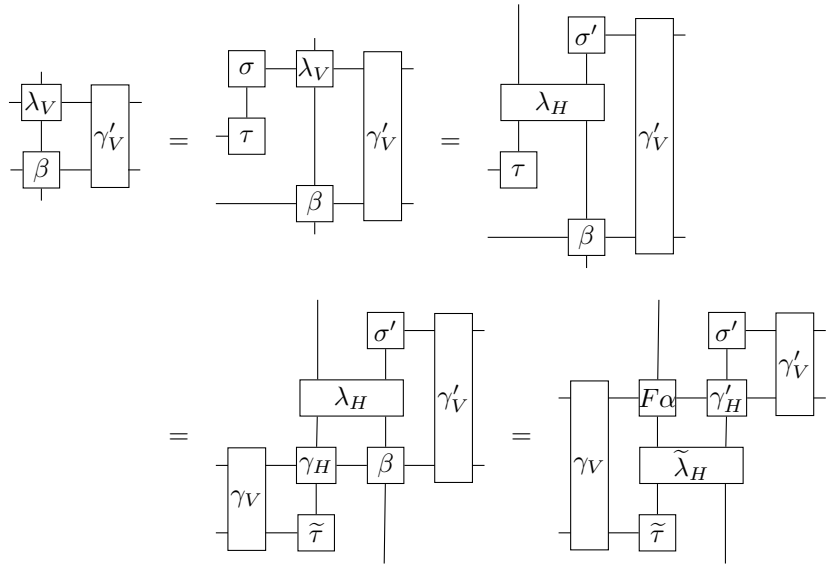
\begin{align*}
\begin{doublestring}
\node[basic] (lambda) {$\lambda_\Ve$};
\node[basic] (beta) [below=of lambda] {$\beta$};
\node[ghost] (gamma'1) [right=of lambda] {};
\node[ghost] (gamma'2) [right=of beta] {};
\node[h2] at ($(gamma'1)!1/2!(gamma'2)$) {$\gamma'_\Ve$};
\inedge{lambda}{beta}{}
\inedge{lambda}{gamma'1}{}
\inedge{beta}{gamma'2}{}
\outabove{lambda}{}
\outright{gamma'1}{}
\outright{gamma'2}{}
\outbelow{beta}{}
\outleft{beta}{}
\outleft{lambda}{}
\end{doublestring}
&
=
\begin{doublestring}
\node[basic] (lambda) {\makebox[0mm]{$\lambda_\Ve$}};
\node[ghost] (gamma'1) [right=of lambda] {};
\node[ghost] (gamma'2) [below=of gamma'1] {};
\node[ghost] (gamma'3) [below=of gamma'2] {};
\node[h3] at (gamma'2) {$\gamma'_\Ve$};
\node[basic] (beta) [left=of gamma'3] {$\beta$};
\node[basic] (sigma) [left=of lambda] {$\sigma$};
\node[basic] (tau) [below=of sigma] {$\tau$};
\node[ghost] (tauout) [below=of tau] {};
\inedge{sigma}{lambda}{}
\inedge{sigma}{tau}{}
\inedge{lambda}{beta}{}
\inedge{lambda}{gamma'1}{}
\inedge{beta}{gamma'3}{}
\outabove{lambda}{}
\outright{gamma'1}{}
\outright{gamma'3}{}
\outbelow{beta}{}
\outleftx{beta}{tauout}{}
\outleft{tau}{}
\end{doublestring}
=
\begin{doublestring}
\node[basic] (sigma') {$\sigma'$};
\node[ghost] (sigma'out) [left=of sigma'] {};
\node[ghost] (gamma'1) [right=of sigma'] {};
\node[ghost] (gamma'2) [below=of gamma'1] {};
\node[ghost] (gamma'3) [below=of gamma'2] {};
\node[ghost] (gamma'4) [below=of gamma'3] {};
\node[h4] at ($(gamma'1)!1/2!(gamma'4)$) {$\gamma'_\Ve$};
\node[basic] (beta) [left=of gamma'4] {$\beta$};
\node[ghost] (lambda2) [below=of sigma'] {};
\node[ghost] (lambda1) [left=of lambda2] {};
\node[w2] at ($(lambda1)!1/2!(lambda2)$) {$\lambda_\Ho$};
\node[basic] (tau) [below=of lambda1] {$\tau$};
\node[ghost] (tauout) [below=of tau] {};
\inedge{sigma'}{gamma'1}{}
\inedge{sigma'}{lambda2}{}
\inedge{lambda1}{tau}{}
\inedge{lambda2}{beta}{}
\inedge{beta}{gamma'4}{}
\outabovex{lambda1}{sigma'out}{}
\outright{gamma'1}{}
\outright{gamma'4}{}
\outbelow{beta}{}
\outleftx{beta}{tauout}{}
\outleft{tau}{}
\end{doublestring}
\\
&
=
\begin{doublestring}
\node[basic] (sigma') {$\sigma'$};
\node[ghost] (sigma'out) [left=of sigma'] {};
\node[ghost] (gamma'1) [right=of sigma'] {};
\node[ghost] (gamma'2) [below=of gamma'1] {};
\node[ghost] (gamma'3) [below=of gamma'2] {};
\node[h3] at (gamma'2) {$\gamma'_\Ve$};
\node[basic] (beta) [left=of gamma'3] {$\beta$};
\node[ghost] (lambda2) [left=of gamma'2] {};
\node[ghost] (lambda1) [left=of lambda2] {};
\node[w2] at ($(lambda1)!1/2!(lambda2)$) {$\lambda_\Ho$};
\node[basic] (gammaH) [left=of beta] {$\gamma_\Ho$};
\node[basic] (taut) [below=of gammaH] {$\twid{\tau}$};
\node[ghost] (tautout) [right=of taut] {};
\node[ghost] (gammaV1) [left=of gammaH] {};
\node[ghost] (gammaV2) [below=of gammaV1] {};
\node[h2] at ($(gammaV1)!1/2!(gammaV2)$) {$\gamma_\Ve$};
\inedge{sigma'}{gamma'1}{}
\inedge{sigma'}{lambda2}{}
\inedge{lambda1}{gammaH}{}
\inedge{lambda2}{beta}{}
\inedge{beta}{gamma'3}{}
\inedge{gammaH}{beta}{}
\inedge{gammaH}{taut}{}
\inedge{gammaV1}{gammaH}{}
\inedge{gammaV2}{taut}{}
\outabovex{lambda1}{sigma'out}{}
\outright{gamma'1}{}
\outright{gamma'3}{}
\outbelowx{beta}{tautout}{}
\outleft{gammaV2}{}
\outleft{gammaV1}{}
\end{doublestring}
=
\begin{doublestring}
\node[basic] (sigma') {$\sigma'$};
\node[ghost] (sigma'out) [left=of sigma'] {};
\node[ghost] (gamma'V1) [right=of sigma'] {};
\node[ghost] (gamma'V2) [below=of gamma'V1] {};
\node[h2] at ($(gamma'V1)!1/2!(gamma'V2)$) {$\gamma'_\Ve$};
\node[basic] (gamma'H) [below=of sigma'] {$\gamma'_\Ho$};
\node[basic] (Falpha) [left=of gamma'H] {\makebox[0mm]{$F\alpha$}};
\node[ghost] (lambdat1) [below=of Falpha] {};
\node[ghost] (lambdat2) [right=of lambdat1] {};
\node[w2] at ($(lambdat1)!1/2!(lambdat2)$) {$\twid{\lambda}_\Ho$};
\node[basic] (taut) [below=of lambdat1] {$\twid{\tau}$};
\node[ghost] (tautout) [right=of taut] {};
\node[ghost] (gammaV1) [left=of Falpha] {};
\node[ghost] (gammaV2) [below=of gammaV1] {};
\node[ghost] (gammaV3) [below=of gammaV2] {};
\node[h3] at (gammaV2) {$\gamma_\Ve$};
\inedge{sigma'}{gamma'V1}{}
\inedge{sigma'}{gamma'H}{}
\inedge{gamma'H}{gamma'V2}{}
\inedge{Falpha}{gamma'H}{}
\inedge{Falpha}{lambdat1}{}
\inedge{gamma'H}{lambdat2}{}
\inedge{lambdat1}{taut}{}
\inedge{gammaV1}{Falpha}{}
\inedge{gammaV3}{taut}{}
\outabovex{Falpha}{sigma'out}{}
\outright{gamma'1}{}
\outright{gamma'2}{}
\outbelowx{lambdat2}{tautout}{}
\outleft{gammaV3}{}
\outleft{gammaV1}{}
\end{doublestring}
\\
&
=
\begin{doublestring}
\node[basic] (sigmat') {\makebox[0mm]{$\twid{\sigma}'$}};
\node[ghost] (sigmat'out) [above=of sigmat'] {};
\node[ghost] (lambdat2) [below=of sigmat'] {};
\node[ghost] (lambdat1) [left=of lambdat2] {};
\node[w2] at ($(lambdat1)!1/2!(lambdat2)$) {$\twid{\lambda}_\Ho$};
\node[basic] (taut) [below=of lambdat1] {$\twid{\tau}$};
\node[ghost] (tautout) [right=of taut] {};
\node[ghost] (gammaV4) [left=of taut] {};
\node[ghost] (gammaV3) [above=of gammaV4] {};
\node[ghost] (gammaV2) [above=of gammaV3] {};
\node[ghost] (gammaV1) [above=of gammaV2] {};
\node[h4] at ($(gammaV1)!1/2!(gammaV4)$) {$\gamma_\Ve$};
\node[basic] (Falpha) [right=of gammaV1] {\makebox[0mm]{$F\alpha$}};
\inedge{gammaV1}{Falpha}{}
\inedge{gammaV4}{taut}{}
\inedge{Falpha}{lambdat1}{}
\inedge{lambdat1}{taut}{}
\inedge{sigmat'}{lambdat2}{}
\outabove{Falpha}{}
\outrightx{Falpha}{sigmat'out}{}
\outright{sigmat'}{}
\outbelowx{lambdat2}{tautout}{}
\outleft{gammaV4}{}
\outleft{gammaV1}{}
\end{doublestring}
=
\begin{doublestring}
\node[basic] (sigmat') {\makebox[0mm]{$\twid{\sigma}'$}};
\node[ghost] (sigmat'out) [above=of sigmat'] {};
\node[basic] (taut') [below=of sigmat'] {\makebox[0mm]{$\twid{\tau}'$}};
\node[basic] (lambdat) [left=of taut']
{\makebox[0mm]{$\twid{\lambda}_\Ve$}};
\node[ghost] (gammaV3) [left=of lambdat] {};
\node[ghost] (gammaV2) [above=of gammaV3] {};
\node[ghost] (gammaV1) [above=of gammaV2] {};
\node[h3] at (gammaV2) {$\gamma_\Ve$};
\node[basic] (Falpha) [right=of gammaV1] {\makebox[0mm]{$F\alpha$}};
\inedge{gammaV1}{Falpha}{}
\inedge{gammaV3}{lambdat}{}
\inedge{Falpha}{lambdat}{}
\inedge{lambdat}{taut'}{}
\inedge{sigmat'}{taut'}{}
\outabove{Falpha}{}
\outrightx{Falpha}{sigmat'out}{}
\outright{sigmat'}{}
\outbelow{lambdat}{}
\outleft{gammaV3}{}
\outleft{gammaV1}{}
\end{doublestring}
\\
&
=
\begin{doublestring}
\node[basic] (Falpha) {$F\alpha$};
\node[basic] (lambdat) [below=of Falpha] {\raisebox{-1mm}{$\twid{\lambda}_\Ve$}};
\node[ghost] (gamma1) [left=of Falpha] {};
\node[ghost] (gamma2) [left=of lambdat] {};
\node[h2] at ($(gamma1)!1/2!(gamma2)$) {$\gamma_\Ve$};
\inedge{Falpha}{lambdat}{}
\inedge{gamma1}{Falpha}{}
\inedges{gamma2}{lambdat}{}
\outabove{Falpha}{}
\outright{Falpha}{}
\outright{lambdat}{}
\outbelow{lambdat}{}
\outleft{gamma2}{}
\outleft{gamma1}{}
\end{doublestring},
\end{align*}
where the first and last equalities are by the companionship
identities~\eqref{eq:companion}, the second and second to last are
by~\eqref{eq:dbl-coh-1hc}, the third and third to last are
by~\eqref{eq:dbl-coh-1vc}, and the middle equality is
by~\eqref{eq:dbl-coh-2a}. 

\bibliography{mathrefs}

\begin{thebibliography}{10}

\bibitem{ArMc}
N.~Arkor and D.~McDermott.
\newblock The formal theory of relative monads.
\newblock {\em Journal of Pure and Applied Algebra}, 228(9):107676, 2024.

\bibitem{BKP}
R.~Blackwell, G.~M. Kelly, and A.~J. Power.
\newblock Two-dimensional monad theory.
\newblock {\em Journal of Pure and Applied Algebra}, 59:1--41, 1989.

\bibitem{BrowAHT}
K.~S. Brown.
\newblock Abstract homotopy theory and generalized sheaf cohomology.
\newblock {\em Transactions of the American Mathematical Society},
  186:419--458, 1973.

\bibitem{Camp}
A.~Campbell.
\newblock The gregarious model structure for double categories.
\newblock Seminar slides, Masaryk University Algebra Seminar, 18 June 2020,
  available at \url{https://acmbl.github.io}, 2020.

\bibitem{Cott}
T.~Cottrell.
\newblock A study of {P}enon weak $n$-categories, part 2: a multisimplicial
  nerve construction.
\newblock {\em Cahiers de Topologie et G{\'e}om{\'e}trie Diff{\'e}rentielle
  Cat{\'e}goriques}, 45(1):32--114, 2019.

\bibitem{DPP}
R.~Dawson, R.~Par{\'e}, and D.~Pronk.
\newblock The span construction.
\newblock {\em Theory and Applications of Categories}, 24(13):302--377, 2010.

\bibitem{GrPaADC}
M.~Grandis and R.~Par{\'e}.
\newblock Adjoint for double categories.
\newblock {\em Cahiers de Topologie et G{\'e}om{\'e}trie Diff{\'e}rentielle
  Cat{\'e}goriques}, 45(3):193--240, 2004.

\bibitem{KeStRE2}
G.~M. Kelly and R.~Street.
\newblock Review of the elements of 2-categories.
\newblock In G.~M. Kelly, editor, {\em Category Seminar (Proceedings Sydney
  Category Theory Seminar 1972/1973)}, pages 75--103. Springer, Berlin, 1974.

\bibitem{LackQMS}
S.~Lack.
\newblock A {Q}uillen model structure for bicategories.
\newblock {\em K-Theory}, 33:185--197, 2004.

\bibitem{SDN}
T.~Leinster.
\newblock A survey of definitions of $n$-category.
\newblock {\em Theory and Applications of Categories}, 10(1):1--70, 2002.

\bibitem{MakkFOL}
M.~Makkai.
\newblock First order logic with dependent sorts, with applications to category
  theory.
\newblock Preprint, available at \url{https://math.mcgill.ca/makkai}, 1995.

\bibitem{MSV}
L.~Moser, M.~Sarazola, and P.~Verdugo.
\newblock Double categorical equivalences.
\newblock In preparation, 2025.

\bibitem{MyerSDD}
D.~J. Myers.
\newblock String diagrams for double categories and equipments.
\newblock Preprint \href{https://arxiv.org/abs/1612.02762}{arXiv:1612.02762},
  available at arXiv.org, 2016.

\bibitem{Pell}
R.~Pellissier.
\newblock {\em Cat{\'e}gories enrichies faibles}.
\newblock PhD thesis, Universit{\'e} de Nice, 2002.

\bibitem{SaraDCE}
M.~Sarazola.
\newblock Double categorical equivalences.
\newblock Seminar slides, Category Theory 2025, Brno, available at
  \url{https://conference.math.muni.cz/ct2025}, 2025.

\bibitem{ShulCCL}
M.~Shulman.
\newblock Comparing composites of left and right derived functors.
\newblock {\em New York Journal of Mathematics}, 17:75--125, 2011.

\bibitem{SimpCMS}
C.~Simpson.
\newblock A closed model structure for $n$-categories, internal $\mathit{Hom}$,
  $n$-stacks and generalized {S}eifert--{V}an {K}ampen.
\newblock Preprint
  \href{https://arxiv.org/abs/alg-geom/9704006}{arXiv:alg-geom/9704006},
  available at arXiv.org, 1997.

\bibitem{SimpHTH}
C.~Simpson.
\newblock {\em Homotopy Theory of Higher Categories}, volume~19 of {\em New
  Mathematical Monographs}.
\newblock Cambridge University Press, Cambridge, 2011.

\bibitem{StreCS}
R.~Street.
\newblock Categorical structures.
\newblock In M.~Hazewinkel, editor, {\em Handbook of Algebra}, volume~1, pages
  529--577. Elsevier, Amsterdam, 1996.

\end{thebibliography}

\end{document}